\documentclass[11pt,journal,onecolumn,draftcls]{IEEEtran}
\usepackage{mathrsfs}
\usepackage{bbm}
\usepackage{threeparttable}
\usepackage[chatter]{rotating}

\IEEEoverridecommandlockouts

\usepackage{bbding}
\usepackage{cite}
\usepackage{graphicx}
\usepackage[cmex10]{amsmath}
\usepackage{amsfonts}
\usepackage{algorithmic}
\usepackage{algorithm}
\usepackage{array}
\usepackage{mdwmath}
\usepackage{mdwtab}
\usepackage{eqparbox}
\usepackage[caption=false,font=footnotesize]{subfig}
\usepackage{fixltx2e}
\usepackage{stfloats}
\usepackage{amssymb}
\usepackage{color}
\usepackage{multirow}
\usepackage{dsfont}
\usepackage{bm}

\newtheorem{theorem}{Theorem}[section]
\newtheorem{lemma}[theorem]{Lemma}

\newtheorem{assumption}[theorem]{Assumption}

\newtheorem{definition}[theorem]{Definition}

\renewcommand{\Re}{\mathfrak{Re}}

\def\mc{\mathcal}
\def\bs{\boldsymbol}
\def\defeq{\triangleq}
\def\st{{\mathrm{subject\;to}}}
\DeclareMathOperator*{\minimize}{\mathrm{minimize}}

\def\Re{\mathrm{Re}}
\def\wt{\widetilde}

\def\E{\mathbb{E}}
\def\col{\mathrm{col}}
\def\diag{\mathrm{diag}}

\def\Tr{\mathrm{Tr}}
\def\vec{\mathrm{vec}}

\def\T{\mathsf{T}}
\def\N{\mathcal{N}}

\begin{document}

\title{Diffusion Adaptation Over Networks Under \\ Imperfect Information Exchange and \\  Non-Stationary Data}
%
\author{Xiaochuan~Zhao,~\IEEEmembership{Student~Member,~IEEE,}
        Sheng-Yuan~Tu,~\IEEEmembership{Student~Member,~IEEE,}
        and Ali~H.~Sayed,~\IEEEmembership{Fellow,~IEEE}
\thanks{The authors are with Department of Electrical Engineering, University of California, Los Angeles, CA 90095,
Email: \{xzhao,~shinetu,~sayed\}@ee.ucla.edu.}
\thanks{This work was supported in part by NSF grants CCF-0942936 and CCF-1011918. A short version of this work was accepted for
the IEEE International Conference on Communications (ICC), Ottawa, ON, Canada, June 2012 \cite{Zhao12ICC}.}}

%


%
\maketitle

\begin{abstract}
Adaptive networks rely on in-network and collaborative processing among distributed agents to deliver enhanced performance in estimation and inference tasks. Information is exchanged among the nodes, usually over noisy links. The combination weights that are used by the nodes to fuse information from their neighbors play a critical role in influencing the adaptation and tracking abilities of the network. This paper first investigates the mean-square performance of general adaptive diffusion algorithms in the presence of various sources of imperfect information exchanges, quantization errors, and model non-stationarities. Among other results, the analysis reveals that link noise over the regression data modifies the dynamics of the network evolution in a distinct way, and leads to biased estimates in steady-state. The analysis also reveals how the network mean-square performance is dependent on the combination weights. We use these observations to show how the combination weights can be optimized and adapted. Simulation results illustrate the theoretical findings and match well with theory.
\end{abstract}


\newpage

\section{Introduction}
\IEEEPARstart{A}{n} adaptive network consists of a collection of agents that are interconnected to each other and solve distributed estimation and inference problems in a collaborative manner. Two useful strategies that enable adaptation and learning over such networks in real-time are the incremental strategy \cite{Tsitsiklis84TAC,Tsitsiklis86TAC,Bertsekas97JOP,Nedic01JOP,Rabbat05JSAC,Lopes07TSP} and the diffusion strategy \cite{Lopes06ASAP,Lopes08TSP,Cattivelli08ASILOMAR,Cattivelli10TSP,Cattivelli10TAC}. Incremental strategies rely on the use of a Hamiltonian cycle, i.e., a cyclic path that covers all nodes in the network, which is generally difficult to enforce since determining a Hamiltonian cycle is an NP-hard problem. In addition, cyclic trajectories are not robust to node or link failure. In comparison, diffusion strategies are scalable, robust, and able to match well the performance of incremental networks. In adaptive diffusion implementations, information is processed locally at the nodes and then diffused in real-time across the network. Diffusion strategies were originally proposed in \cite{Lopes06ASAP,Lopes08TSP,Cattivelli08ASILOMAR} and further extended and studied in \cite{Cattivelli08TSP,Li09SSP,Cattivelli10TSP,Cattivelli10TAC,Takahashi10TSP,Takahashi10ICASSP,Chouvardas11TSP}. They have been applied to model self-organized and complex behavior encountered in biological networks, such as fish schooling \cite{Tu11JSTSP}, bird flight formations \cite{Cattivelli11TSP}, and bee swarming \cite{Li11EURASIP}. Diffusion strategies have also been applied to online learning of Gaussian mixture models \cite{Towfic11MLSP,Weng11Sensors} and to general distributed optimization problems \cite{Chen11TSP}. There have also been several useful works in the literature on distributed consensus-type strategies, with application to multi-agent formations and distributed processing \cite{Jadbabaie03TAC,Fax04TAC,Olfati06TAC,Olfati04TAC,Barbarossa07SPM,Nedic09TSP,Dimakis10PROC,Kar11TSP}. The main difference between these works and the diffusion approach of \cite{Lopes08TSP,Cattivelli10TSP,Cattivelli10TAC} is the latter's emphasis on the role of adaptation and learning over networks.

In the original diffusion least-mean-squares (LMS) strategy \cite{Lopes08TSP,Cattivelli10TSP}, the weight estimates that are exchanged among the nodes can be subject to quantization errors and additive noise over the communication links. Studying the degradation in mean-square performance that results from these particular perturbations can be pursued, for both incremental and diffusion strategies, by extending the mean-square analysis already presented in \cite{Lopes08TSP,Cattivelli10TSP}, in the same manner that the tracking analysis of conventional stand-alone adaptive filters was obtained from the counterpart results in the stationary case (as explained in
\cite[Ch. 21]{Sayed08}). Useful results along these lines, which study the effect of link noise during the exchange of the weight estimates, already appear for the traditional diffusion algorithm in the works \cite{Abdolee11DCOSS,Khalili11ACSP,Khalili12TSP,Tu11GlobeCom} and for consensus-based algorithms in \cite{Kar09TSP,Mateos09EUROSIP}. In this paper, our objective is to go beyond these earlier studies by taking into account additional effects, and by considering a more general algorithmic structure. The reason for this level of generality is because the analytical results will help reveal which noise sources influence the network performance more seriously, in what manner, and at what stage of the adaptation process. The results will suggest important remedies and mechanisms to adapt the combination weights in real-time. Some of these insights are hard to get if one focuses solely on noise during the exchange of the weight estimates. The analysis will further show that noise during the exchange of the regression data plays a more critical role than other sources of imperfection: this particular noise alters the learning dynamics and modes of the network, and biases the weight estimates. Noises related to the exchange of other pieces of information do not alter the dynamics of the network but contribute to the deterioration of the network performance.

To arrive at these results, in this paper, we first consider a generalized analysis that applies to a broad class of diffusion adaptation strategies (see \eqref{eqn:idealdiffusionpriordiff}--\eqref{eqn:idealdiffusionpostdiff} further ahead; this class includes the original diffusion strategies \eqref{eqn:ctaalgorithm} and \eqref{eqn:atcalgorithm} as two special cases). The analysis allows us to account for various sources of information noise over the communication links. We allow for noisy exchanges during \emph{each} of the three processing steps of the adaptive diffusion algorithm (the two combination steps \eqref{eqn:idealdiffusionpriordiff} and \eqref{eqn:idealdiffusionpostdiff} and the adaptation step \eqref{eqn:idealdiffusionincremental}). In this way, we are able to examine how the three sets of combination coefficients $\{a_{1,lk},c_{lk},a_{2,lk}\}$ in \eqref{eqn:idealdiffusionpriordiff}--\eqref{eqn:idealdiffusionpostdiff} influence the propagation of the noise signals through the network dynamics. Our results further reveal how the network mean-square-error performance is dependent on these combination weights. Following this line of reasoning, the analysis leads to algorithms \eqref{eqn:relativevariance} and \eqref{eqn:relativevarianceadapt} further ahead for choosing the combination coefficients to improve the steady-state network performance.

It should be noted that several combination rules, such as the Metropolis rule \cite{Metropolis53JCP} and the maximum degree rule \cite{Xiao05IPSN}, were proposed previously in the literature --- especially in the context of consensus-based iterations \cite{Xiao04SCL,Xiao05IPSN,Jakovetic10TSP}. These schemes, however, usually suffer performance degradation in the presence of noisy information exchange since they ignore the network noise profile \cite{Takahashi10TSP}. When the noise variance differs across the nodes, it becomes necessary to design combination rules that are aware of this variation as outlined further ahead in Section VI-B. Moreover, in a mobile network \cite{Tu11JSTSP} where nodes are on the move and where neighborhoods evolve over time, it is even more critical to employ adaptive combination strategies that are able to track the variations in the noise profile in order to cope with such dynamic environments. This issue is taken up in Section VI-C.

\subsection{Notation}
We use lowercase letters to denote vectors, uppercase letters for matrices, plain letters for deterministic variables, and boldface letters for random variables. We also use $(\cdot)^*$ to denote conjugate transposition, ${\mathrm{Tr}}(\cdot)$ for the trace of its matrix argument, $\rho(\cdot)$ for the spectral radius of its matrix argument, $\otimes$ for the Kronecker product, and ${\mathrm{vec}}(\cdot)$ for a vector formed by stacking the columns of its matrix argument. We further use $\diag\{\cdots\}$ to denote a (block) diagonal matrix formed from its arguments, and $\col\{\cdots\}$ to denote a column vector formed by stacking its arguments on top of each other. All vectors in our treatment are column vectors, with the exception of the regression vectors, ${\bs{u}}_{k,i}$, and the associated noise signals, ${\bs{v}}_{lk,i}^{(u)}$, which are taken to be row vectors for convenience of presentation.

\section{Diffusion Algorithms with Imperfect Information Exchange}
We consider a connected network consisting of $N$ nodes. Each node $k$ collects scalar measurements ${\bs{d}}_k(i)$ and $1\times M$ regression data vectors ${\bs{u}}_{k,i}$ over successive time instants $i\ge0$. Note that we use parenthesis to refer to the time-dependence of scalar variables, as in ${\bs{d}}_k(i)$, and subscripts to refer to the time-dependence of vector variables, as in ${\bs{u}}_{k,i}$. The measurements across all nodes are assumed to be related to an unknown $M\times1$ vector $w^o$ via a linear regression model of the form \cite{Sayed08}:
\begin{align}
\label{eqn:lineardatamodel}
{\bs{d}}_{k}(i)={\bs{u}}_{k,i}w^o+{\bs{v}}_k(i)
\end{align}
where ${\bs{v}}_k(i)$ denotes the measurement or model noise with zero mean and variance $\sigma_{v,k}^{2}$. The vector $w^o$ in \eqref{eqn:lineardatamodel} denotes the parameter of interest, such as the parameters of some underlying physical phenomenon, the taps of a communication channel, or the location of food sources or predators. Such data models are also useful in studies on hybrid combinations of adaptive filters \cite{Arenas05TIM,Mandic07ICASSP,Silva08TSP,Kozat09ICASSP,Candido10TSP}.

The nodes in the network would like to estimate $w^o$ by solving the following minimization problem:
\begin{align}
\label{eqn:globaloptimization}
\minimize_{w} \quad \sum_{k=1}^{N} {\mathbb{E}}|{\bs{d}}_k(i)-{\bs{u}}_{k,i}w|^2
\end{align}
In previous works \cite{Lopes08TSP,Cattivelli08TSP,Cattivelli10TSP}, we introduced and studied several distributed strategies of the diffusion type that allow nodes to cooperate with each other in order to solve problems of the form \eqref{eqn:globaloptimization} in an adaptive manner. These diffusion strategies endow networks with adaptation and learning abilities, and enable information to diffuse through the network in real-time. We review the adaptive diffusion strategies below.

\subsection{Diffusion Adaptation with Perfect Information Exchange}
In \cite{Lopes08TSP, Cattivelli10TSP}, two classes of diffusion algorithms were proposed. One class is the so-called Combine-then-Adapt (CTA) strategy:
\begin{align}
\label{eqn:ctaalgorithm}
\left\{\begin{aligned}
{\bs{\phi}}_{k,i-1}&=\sum_{l\in{\mc{N}}_k}a_{1,lk}{\bs{w}}_{l,i-1}\\
{\bs{w}}_{k,i}&={\bs{\phi}}_{k,i-1}\!+\!\mu_k\sum_{l\in{\mc{N}}_k}c_{lk}{\bs{u}}_{l,i}^*
[{\bs{d}}_{l}(i)\!-\!{\bs{u}}_{l,i}{\bs{\phi}}_{k,i-1}]\\
\end{aligned}\right.
\end{align}
and the second class is the so-called Adapt-then-Combine (ATC) strategy:
\begin{align}
\label{eqn:atcalgorithm}
\left\{\begin{aligned}
{\bs{\psi}}_{k,i}&={\bs{w}}_{k,i-1}\!+\!\mu_k\sum_{l\in{\mc{N}}_k}c_{lk}{\bs{u}}_{l,i}^*
[{\bs{d}}_{l}(i)\!-\!{\bs{u}}_{l,i}{\bs{w}}_{k,i-1}]\\
{\bs{w}}_{k,i}&=\sum_{l\in{\mc{N}}_k}a_{2,lk}{\bs{\psi}}_{l,i}\\
\end{aligned}\right.
\end{align}
where the $\{\mu_k\}$ are small positive step-size parameters and the $\{a_{1,lk},c_{lk},a_{2,lk}\}$ are nonnegative entries of the $N\times N$ matrices $\{A_1,C,A_2\}$, respectively. The coefficients $\{a_{1,lk},c_{lk},a_{2,lk}\}$ are zero whenever node $l$ is not connected to node $k$, i.e., $l\notin{\mc{N}}_k$, where ${\mc{N}}_k$ denotes the neighborhood of node $k$. The two strategies \eqref{eqn:ctaalgorithm} and \eqref{eqn:atcalgorithm} can be integrated into one broad class of diffusion adaptation \cite{Cattivelli10TSP}:
\begin{align}
\label{eqn:idealdiffusionpriordiff}
{\bs{\phi}}_{k,i-1}&=\sum_{l\in{\mc{N}}_k}a_{1,lk}{\bs{w}}_{l,i-1}\\
\label{eqn:idealdiffusionincremental}
{\bs{\psi}}_{k,i}&={\bs{\phi}}_{k,i-1}+\mu_k\sum_{l\in{\mc{N}}_k}c_{lk}{\bs{u}}_{l,i}^*[{\bs{d}}_{l}(i)-{\bs{u}}_{l,i}{\bs{\phi}}_{k,i-1}]\\
\label{eqn:idealdiffusionpostdiff}
{\bs{w}}_{k,i}&=\sum_{l\in{\mc{N}}_k}a_{2,lk}{\bs{\psi}}_{l,i}
\end{align}
Several diffusion strategies can be obtained as special cases of \eqref{eqn:idealdiffusionpriordiff}--\eqref{eqn:idealdiffusionpostdiff} through proper selection of the coefficients $\{a_{1,lk},c_{lk},a_{2,lk}\}$. For example, to recover the CTA strategy \eqref{eqn:ctaalgorithm}, we set $A_2=I_N$, and to recover the ATC strategy \eqref{eqn:atcalgorithm}, we set $A_1=I_N$, where $I_N$ denotes the $N\times N$ identity matrix. In the general diffusion strategy \eqref{eqn:idealdiffusionpriordiff}--\eqref{eqn:idealdiffusionpostdiff}, each node $k$ evaluates its estimate ${\bs{w}}_{k,i}$ at time $i$ by relying solely on the data collected from its neighbors through steps \eqref{eqn:idealdiffusionpriordiff} and \eqref{eqn:idealdiffusionpostdiff} and on its local measurements through step \eqref{eqn:idealdiffusionincremental}. The matrices $A_1$, $A_2$, and $C$ are required to be left or right-stochastic, i.e.,
\begin{align}
\label{eqn:A1A2def}
A_1^{\T}{\mathds{1}}_N={\mathds{1}}_N,\qquad A_2^{\T}{\mathds{1}}_N={\mathds{1}}_N,\qquad
C{\mathds{1}}_N={\mathds{1}}_N
\end{align}
where ${\mathds{1}}_N$ denotes the $N\times1$ vector whose entries are all one. This means that each node performs a convex combination of the estimates received from its neighbors at every iteration $i$.

The mean-square performance and convergence properties of the diffusion algorithm \eqref{eqn:idealdiffusionpriordiff}--\eqref{eqn:idealdiffusionpostdiff} have already been studied in detail in \cite{Lopes08TSP, Cattivelli10TSP}. For the benefit of the analysis in the subsequent sections, we present below in \eqref{eqn:idealerrorrecursion} the recursion describing the evolution of the weight error vectors across the network. To do so, we introduce the error vectors:
\begin{align}
\label{eqn:tildephikidef}
{\widetilde{\bs{\phi}}}_{k,i-1}&\defeq w^o-{\bs{\phi}}_{k,i-1}\\
\label{eqn:tildepsikidef}
{\widetilde{\bs{\psi}}}_{k,i}&\defeq w^o-{\bs{\psi}}_{k,i}\\
\label{eqn:tildewkidef}
{\widetilde{\bs{w}}}_{k,i}&\defeq w^o-{\bs{w}}_{k,i}
\end{align}
and substitute the linear model \eqref{eqn:lineardatamodel} into the adaptation step \eqref{eqn:idealdiffusionincremental} to find that
\begin{align}
\label{eqn:idealdiffusionincremental2}
{\wt{\bs{\psi}}}_{k,i}=(I_M-\mu_k{\bs{R}}_{k,i}){\wt{\bs{\phi}}}_{k,i-1}-\mu_k\sum_{l\in\N_k}
c_{lk}{\bs{s}}_{l,i}
\end{align}
where the $M\times M$ matrix ${\bs{R}}_{k,i}$ and the $M\times1$ vector ${\bs{s}}_{k,i}$ are defined as:
\begin{align}
\label{eqn:rndRkidef}
{\bs{R}}_{k,i}&\defeq\sum_{l\in\N_k}c_{lk}{\bs{u}}_{l,i}^*{\bs{u}}_{l,i}\\
\label{eqn:skidef}
{\bs{s}}_{k,i}&\defeq{\bs{u}}_{k,i}^*{\bs{v}}_k(i)
\end{align}
We further collect the various quantities across all nodes in the network into the following block vectors and matrices:
\begin{align}
\label{eqn:rndbigRidef}
{\bs{\mc{R}}}_i&\defeq{\mathrm{diag}}\left\{{\bs{R}}_{1,i},\dots,{\bs{R}}_{N,i}\right\}\\
\label{eqn:sidef}
{\bs{s}}_i&\defeq{\mathrm{col}}\left\{{\bs{s}}_{1,i},\dots,{\bs{s}}_{N,i}\right\}\\
\label{eqn:bigMdef}
{\mc{M}}&\defeq{\mathrm{diag}}\left\{\mu_1I_M,\dots,\mu_NI_M\right\}\\
\label{eqn:phidef}
{\wt{\bs{\phi}}}_i&\defeq{\mathrm{col}}\left\{{\wt{\bs{\phi}}}_{1,i},\dots,{\wt{\bs{\phi}}}_{N,i}\right\}\\
\label{eqn:psidef}
{\wt{\bs{\psi}}}_i&\defeq{\mathrm{col}}\left\{{\wt{\bs{\psi}}}_{1,i},\dots,{\wt{\bs{\psi}}}_{N,i}\right\}\\
\label{eqn:widef}
{\wt{\bs{w}}}_i&\defeq{\mathrm{col}}\left\{{\wt{\bs{w}}}_{1,i},\dots,{\wt{\bs{w}}}_{N,i}\right\}
\end{align}
Then, from \eqref{eqn:idealdiffusionpriordiff}, \eqref{eqn:idealdiffusionpostdiff}, and \eqref{eqn:idealdiffusionincremental2}, the recursion for the network error vector ${\wt{\bs{w}}}_i$ is given by
\begin{equation}
\boxed{
\label{eqn:idealerrorrecursion}
{\widetilde{\bs{w}}}_i={\mc{A}}_2^{\T}(I_{NM}-{\mc{M}}{\bs{\mc{R}}}_i){\mc{A}}_1^{\T}{\widetilde{\bs{w}}}_{i-1}
-{\mc{A}}_2^{\T}{\mc{M}}{\mc{C}}^{\T}{\bs{s}}_i
}\end{equation}
where
\begin{align}
\label{eqn:bigA1CA2def}
{\mc{A}}_1\defeq A_1\otimes I_M,\quad
{\mc{C}}\defeq C\otimes I_M,\quad
{\mc{A}}_2\defeq A_2\otimes I_M
\end{align}

\subsection{Noisy Information Exchange}
Each of the steps in \eqref{eqn:idealdiffusionpriordiff}--\eqref{eqn:idealdiffusionpostdiff} involves the sharing of information between node $k$ and its neighbors. For example, in the first step \eqref{eqn:idealdiffusionpriordiff}, all neighbors of node $k$ send their estimates ${\bs{w}}_{l,i-1}$ to node $k$. This transmission is generally subject to additive noise and possibly quantization errors. Likewise, steps \eqref{eqn:idealdiffusionincremental} and \eqref{eqn:idealdiffusionpostdiff} involve the sharing of other pieces of information with node $k$. These exchange steps can all be subject to perturbations (such as additive noise and quantization errors). One of the objectives of this work is to analyze the \emph{aggregate} effect of these perturbations on general diffusion strategies of the type \eqref{eqn:idealdiffusionpriordiff}--\eqref{eqn:idealdiffusionpostdiff} and to propose choices for the combination weights in order to enhance the mean-square performance of the network in the presence of these disturbances.

\begin{figure}[t]
\centering
\includegraphics[scale=0.4]{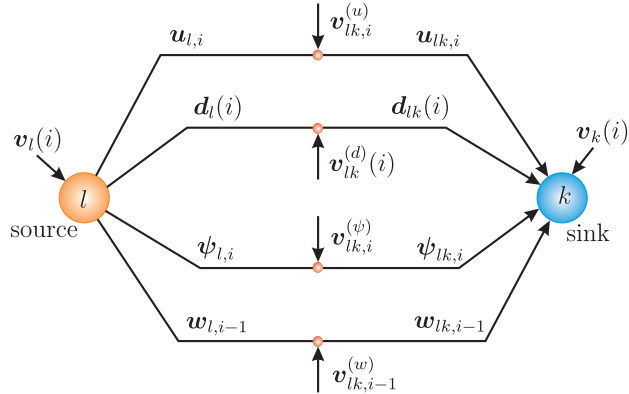}
\caption{Several additive noise sources perturb the exchange of information from node $l$ to node $k$.}
\label{fig:noise}
\end{figure}

So let us examine what happens when information is exchanged over links with additive noise. We model the data received by node $k$ from its neighbor $l$ as
\begin{align}
\label{eqn:noisywlkidef}
{\bs{w}}_{lk,i-1}&\defeq{\bs{w}}_{l,i-1}+{\bs{v}}_{lk,i-1}^{(w)}\\
\label{eqn:noisypsilkidef}
{\bs{\psi}}_{lk,i}&\defeq{\bs{\psi}}_{l,i}+{\bs{v}}_{lk,i}^{(\psi)}\\
\label{eqn:noisydlkidef}
{\bs{d}}_{lk}(i)&\defeq{\bs{d}}_{l}(i)+{\bs{v}}_{lk}^{(d)}(i)\\
\label{eqn:noisyulkidef}
{\bs{u}}_{lk,i}&\defeq{\bs{u}}_{l,i}+{\bs{v}}_{lk,i}^{(u)}
\end{align}
where ${\bs{v}}_{lk,i-1}^{(w)}$ and ${\bs{v}}_{lk,i}^{(\psi)}$ are $M\times1$ noise signals, ${\bs{v}}_{lk,i}^{(u)}$ is a $1\times M$ noise signal, and ${\bs{v}}_{lk}^{(d)}(i)$ is a scalar noise signal (see Fig. \ref{fig:noise}). Observe further that in \eqref{eqn:noisywlkidef}--\eqref{eqn:noisyulkidef}, we are including several sources of information exchange noise. In comparison, references \cite{Abdolee11DCOSS, Khalili11ACSP, Khalili12TSP} only considered the noise source ${\bs{v}}_{lk,i-1}^{(w)}$ in \eqref{eqn:noisywlkidef} and one set of combination coefficients $\{a_{1,lk}\}$; the other coefficients were set to $c_{lk}=a_{2,lk}=0$ for $l\neq k$ and $c_{kk}=a_{2,kk}=1$. In other words, these references only considered \eqref{eqn:noisywlkidef} and the following traditional CTA strategy without exchange of the data $\{{\bs{d}}_l(i),{\bs{u}}_{l,i}\}$ --- compare with \eqref{eqn:ctaalgorithm}; note that the second step in \eqref{eqn:ctaalgorithmsimple} only uses $\{{\bs{d}}_k(i),{\bs{u}}_{k,i}\}$:
\begin{align}
\label{eqn:ctaalgorithmsimple}
\left\{\begin{aligned}
{\bs{\phi}}_{k,i-1}&=\sum_{l\in{\mc{N}}_k}a_{1,lk}{\bs{w}}_{l,i-1}\\
{\bs{w}}_{k,i}&={\bs{\phi}}_{k,i-1}+\mu_k{\bs{u}}_{k,i}^*[{\bs{d}}_{k}(i)-{\bs{u}}_{k,i}{\bs{\phi}}_{k,i-1}]\\
\end{aligned}\right.
\end{align}
The analysis that follows examines the aggregate effect of all four noise sources appearing in \eqref{eqn:noisywlkidef}--\eqref{eqn:noisyulkidef}, in addition to the three sets of combination coefficients appearing in \eqref{eqn:idealdiffusionpriordiff}--\eqref{eqn:idealdiffusionpostdiff}. We introduce the following assumption on the statistical properties of the measurement data and noise signals.
\begin{assumption}[Statistical properties of the variables]
\label{asm:all} \;\; \hskip 1in
\begin{enumerate}
  \item The regression data ${\bs{u}}_{k,i}$ are temporally white and spatially independent random variables with zero mean and covariance matrix $R_{u,k}\defeq\E\,{\bm{u}}_{k,i}^*{\bm{u}}_{k,i}\ge0$.
  \item The noise signals ${\bs{v}}_{k}(i)$, ${\bs{v}}_{lk,i-1}^{(w)}$, ${\bs{v}}_{lk}^{(d)}(i)$, ${\bs{v}}_{lk,i}^{(u)}$, and ${\bs{v}}_{lk,i}^{(\psi)}$ are temporally white and spatially independent random variables with zero mean and (co)variances $\sigma_{v,k}^{2}$, $R_{v,lk}^{(w)}$, $\sigma_{v,lk}^{2}$, $R_{v,lk}^{(u)}$, and $R_{v,lk}^{(\psi)}$, respectively. In addition, $R_{v,lk}^{(w)}$, $\sigma_{v,lk}^{2}$, $R_{v,lk}^{(u)}$, and $R_{v,lk}^{(\psi)}$ are all zero if $l\notin{\mc{N}}_k$ or $l=k$.
  \item The regression data $\{{\bs{u}}_{m,i_1}\}$, the model noise signals $\{{\bs{v}}_{n}(i_2)\}$, and the link noise signals $\{{\bs{v}}_{l_1k_1,j_1}^{(w)}\}$, $\{{\bs{v}}_{l_2k_2}^{(d)}(j_2)\}$, $\{{\bs{v}}_{l_3k_3,j_3}^{(u)}\}$, and $\{{\bs{v}}_{l_4k_4,j_4}^{(\psi)}\}$ are mutually-independent random variables for all $\{i_1,i_2,j_1,j_2,j_3,j_4\}$ and $\{m,n,l_1,l_2,l_3,l_4,k_1,k_2,k_3,k_4\}$.  \hfill \IEEEQED
\end{enumerate}
\end{assumption}
\vspace{0.5\baselineskip}
Using the perturbed data \eqref{eqn:noisywlkidef}--\eqref{eqn:noisyulkidef}, the diffusion algorithm \eqref{eqn:idealdiffusionpriordiff}--\eqref{eqn:idealdiffusionpostdiff} becomes
\begin{align}
\label{eqn:noisydiffusionpriordiffold}
{\bs{\phi}}_{k,i-1}&\!=\!\sum_{l\in{\mc{N}}_k}a_{1,lk}{\bs{w}}_{lk,i-1}\\
\label{eqn:noisydiffusionincrementalold}
{\bs{\psi}}_{k,i}&\!=\!{\bs{\phi}}_{k,i-1}\!+\!\mu_k\sum_{l\in{\mc{N}}_k}\!c_{lk}{\bs{u}}_{lk,i}^*
[{\bs{d}}_{lk}(i)\!-\!{\bs{u}}_{lk,i}{\bs{\phi}}_{k,i-1}]\\
\label{eqn:noisydiffusionpostdiffold}
{\bs{w}}_{k,i}&\!=\!\sum_{l\in{\mc{N}}_k}a_{2,lk}{\bs{\psi}}_{lk,i}
\end{align}
where we continue to use the symbols $\{{\bs{\phi}}_{k,i-1},{\bs{\psi}}_{k,i},{\bs{w}}_{k,i}\}$ to avoid an explosion of notation. From \eqref{eqn:noisywlkidef} and \eqref{eqn:noisypsilkidef}, expressions \eqref{eqn:noisydiffusionpriordiffold}--\eqref{eqn:noisydiffusionpostdiffold} can be rewritten as
\begin{align}
\label{eqn:noisydiffusionpriordiff}
{\bs{\phi}}_{k,i-1}&\!=\!\sum_{l\in{\mc{N}}_k}a_{1,lk}{\bs{w}}_{l,i-1}\!+\!{\bs{v}}_{k,i-1}^{(w)}\\
\label{eqn:noisydiffusionincremental}
{\bs{\psi}}_{k,i}&\!=\!{\bs{\phi}}_{k,i-1}\!+\!\mu_k\sum_{l\in{\mc{N}}_k}c_{lk}{\bs{u}}_{lk,i}^*
[{\bs{d}}_{lk}(i)\!-\!{\bs{u}}_{lk,i}{\bs{\phi}}_{k,i-1}]\\
\label{eqn:noisydiffusionpostdiff}
{\bs{w}}_{k,i}&\!=\!\sum_{l\in{\mc{N}}_k}a_{2,lk}{\bs{\psi}}_{l,i}\!+\!{\bs{v}}_{k,i}^{(\psi)}
\end{align}
where we are introducing the symbols ${\bs{v}}_{k,i-1}^{(w)}$ and ${\bs{v}}_{k,i}^{(\psi)}$ to denote the aggregate $M\times1$ zero-mean noise signals defined over the neighborhood of node $k$:
\begin{align}
\label{eqn:vkiwdef}
{\bs{v}}_{k,i-1}^{(w)}&\defeq\sum_{l\in{\mc{N}}_k\backslash\{k\}}a_{1,lk}{\bs{v}}_{lk,i-1}^{(w)}\\
\label{eqn:vkipsidef}
{\bs{v}}_{k,i}^{(\psi)}&\defeq\sum_{l\in{\mc{N}}_k\backslash\{k\}}a_{2,lk}{\bs{v}}_{lk,i}^{(\psi)}
\end{align}
with covariance matrices
\begin{align}
\label{eqn:Rvkwdef}
R_{v,k}^{(w)}&\defeq\sum_{l\in\N_k\backslash\{k\}}a_{1,lk}^2R_{v,lk}^{(w)}\\
\label{eqn:Rvkpsidef}
R_{v,k}^{(\psi)}&\defeq\sum_{l\in\N_k\backslash\{k\}}a_{2,lk}^2R_{v,lk}^{(\psi)}
\end{align}
It is worth noting that $R_{v,k}^{(w)}$ and $R_{v,k}^{(\psi)}$ depend on the combination coefficients $\{a_{1,lk}\}$ and $\{a_{2,lk}\}$, respectively. This property will be taken into account when optimizing over $\{a_{1,lk}\}$ and $\{a_{2,lk}\}$ in a later section. We further introduce the following scalar zero-mean noise signal:
\begin{align}
\label{eqn:vlkidef}
{\bs{v}}_{lk}(i)\defeq{\bs{v}}_l(i)+{\bs{v}}_{lk}^{(d)}(i)-{\bs{v}}_{lk,i}^{(u)}w^o
\end{align}
for $l\in\N_k\backslash\{k\}$, whose variance is
\begin{align}
\label{eqn:sigmalkdef}
\sigma_{lk}^2\defeq\sigma_{v,l}^2+\sigma_{v,lk}^2+w^{o*}R_{v,lk}^{(u)}w^o
\end{align}
To unify the notation, we define ${\bs{v}}_{kk}(i)\defeq{\bs{v}}_{k}(i)$. Then, from \eqref{eqn:lineardatamodel}, \eqref{eqn:noisydlkidef}, and \eqref{eqn:noisyulkidef}, it is easy to verify that the noisy data $\{{\bs{d}}_{lk}(i),{\bs{u}}_{lk,i}\}$ are related via
\begin{align}
\label{eqn:noisylineardatamodel}
{\bs{d}}_{lk}(i)={\bs{u}}_{lk,i}w^o+{\bs{v}}_{lk}(i)
\end{align}
for $l\in\N_k$. Continuing with the adaptation step \eqref{eqn:noisydiffusionincremental} and substituting \eqref{eqn:noisylineardatamodel}, we get
\begin{align}
\label{eqn:noisyadaptation}
{\bs{\psi}}_{k,i}&={\bs{\phi}}_{k,i-1}+\mu_k\sum_{l\in{\mc{N}}_k}c_{lk}{\bs{u}}_{lk,i}^*
[{\bs{u}}_{lk,i}{\wt{\bs{\phi}}}_{k,i-1}+{\bs{v}}_{lk}(i)]
\end{align}
Then, we can derive the following error recursion for node $k$ (compare with \eqref{eqn:idealdiffusionincremental2}):
\begin{align}
\label{eqn:noisyerrorupdatepsiphi}
{\widetilde{\bs{\psi}}}_{k,i}&=\left(I_M-\mu_k{\bs{R}}_{k,i}'\right){\widetilde{\bs{\phi}}}_{k,i-1}-\mu_k{\bs{z}}_{k,i}
\end{align}
where the $M\times M$ matrix ${\bs{R}}_{k,i}'$ and the $M\times 1$ vector ${\bs{z}}_{k,i}$ are defined as (compare with \eqref{eqn:rndRkidef} and \eqref{eqn:skidef}):
\begin{align}
\label{eqn:rndRkiprimedef}
{\bs{R}}_{k,i}'&\defeq\sum_{l\in{\mc{N}}_k}c_{lk}{\bs{u}}_{lk,i}^*{\bs{u}}_{lk,i}\\
\label{eqn:zkidef}
{\bs{z}}_{k,i}&\defeq\sum_{l\in{\mc{N}}_k}c_{lk}{\bs{u}}_{lk,i}^*{\bs{v}}_{lk}(i)
\end{align}
We further introduce the block vectors and matrices:
\begin{align}
\label{eqn:rndbigRiprimedef}
{\bs{\mc{R}}}_i'&\defeq{\mathrm{diag}}\left\{{\bs{R}}_{1,i}',\dots,{\bs{R}}_{N,i}'\right\}\\
\label{eqn:zidef}
{\bs{z}}_i&\defeq{\mathrm{col}}\left\{{\bs{z}}_{1,i},\dots,{\bs{z}}_{N,i}\right\}\\
\label{eqn:viwdef}
{\bs{v}}_i^{(w)}&\defeq{\mathrm{col}}\left\{{\bs{v}}_{1,i}^{(w)},\dots,{\bs{v}}_{N,i}^{(w)}\right\}\\
\label{eqn:vipsidef}
{\bs{v}}_i^{(\psi)}&\defeq{\mathrm{col}}\left\{{\bs{v}}_{1,i}^{(\psi)},\dots,{\bs{v}}_{N,i}^{(\psi)}\right\}
\end{align}
and the corresponding covariance matrices for ${\bs{v}}_i^{(w)}$ and ${\bs{v}}_i^{(\psi)}$:
\begin{align}
\label{eqn:bigRvwdef}
{\mc{R}}_v^{(w)}&\defeq{\mathrm{diag}}\left\{R_{v,1}^{(w)},\dots,R_{v,N}^{(w)}\right\}\\
\label{eqn:bigRvpsidef}
{\mc{R}}_v^{(\psi)}&\defeq{\mathrm{diag}}\left\{R_{v,1}^{(\psi)},\dots,R_{v,N}^{(\psi)}\right\}
\end{align}
then, from \eqref{eqn:noisydiffusionpriordiff}, \eqref{eqn:noisydiffusionpostdiff}, and \eqref{eqn:noisyerrorupdatepsiphi}, we arrive at the following recursion for the network weight error vector in the presence of noisy information exchange:
\begin{align}
\label{eqn:noisyerrorrecursion}
{\wt{\bs{w}}}_i&={\mc{A}}_2^{\T}{\wt{\bs{\psi}}}_i\!-\!{\bs{v}}_i^{(\psi)}\nonumber\\
{}&={\mc{A}}_2^{\T}\left[(I_{NM}\!-\!{\mc{M}}{\bs{\mc{R}}}_i'){\wt{\bs{\phi}}}_{i-1}
\!-\!{\mc{M}}{\bs{z}}_i\right]\!-\!{\bs{v}}_i^{(\psi)}\nonumber\\
{}&={\mc{A}}_2^{\T}\left[(I_{NM}\!-\!{\mc{M}}{\bs{\mc{R}}}_i')({\mc{A}}_1^{\T}{\wt{\bs{w}}}_{i-1}\!-\!{\bs{v}}_{i-1}^{(w)})
\!-\!{\mc{M}}{\bs{z}}_i\right]\!-\!{\bs{v}}_i^{(\psi)}
\end{align}
That is,
\begin{equation}
\boxed{
\label{eqn:noisyerrorrecursion1}
\begin{aligned}
{\wt{\bs{w}}}_i={\mc{A}}_2^{\T}(I_{NM}\!-\!{\mc{M}}{\bs{\mc{R}}}_i'){\mc{A}}_1^{\T}{\wt{\bs{w}}}_{i-1}
-{\mc{A}}_2^{\T}(I_{NM}\!-\!{\mc{M}}{\bs{\mc{R}}}_i'){\bs{v}}_{i-1}^{(w)}\!-\!{\mc{A}}_2^{\T}{\mc{M}}{\bs{z}}_i\!-\!{\bs{v}}_i^{(\psi)}
\end{aligned}
}\end{equation}
Compared to the previous error recursion \eqref{eqn:idealerrorrecursion}, the noise terms in \eqref{eqn:noisyerrorrecursion1} consist of three parts:
\begin{itemize}
  \item ${\mc{A}}_2^{\T}\left(I_{NM}-{\mc{M}}{\bs{\mc{R}}}_i'\right){\bs{v}}_{i-1}^{(w)}$ is contributed by the noise introduced at the information-exchange step \eqref{eqn:noisydiffusionpriordiffold} \textit{before} adaptation.
  \item ${\mc{A}}_2^{\T}{\mc{M}}{\bs{z}}_i$ is contributed by the noise introduced at the adaptation step \eqref{eqn:noisydiffusionincrementalold}.
  \item ${\bs{v}}_i^{(\psi)}$ is contributed by the noise introduced at the information-exchange step \eqref{eqn:noisydiffusionpostdiffold} \textit{after} adaptation.
\end{itemize}

\section{Convergence in the Mean with a Bias}
Given the weight error recursion \eqref{eqn:noisyerrorrecursion1}, we are now ready to study the mean convergence condition for the diffusion strategy \eqref{eqn:noisydiffusionpriordiffold}--\eqref{eqn:noisydiffusionpostdiffold} in the presence of disturbances during information exchange under Assumption \ref{asm:all}. Taking expectations of both sides of \eqref{eqn:noisyerrorrecursion1}, we get
\begin{align}
\label{eqn:meanrecursion}
{\mathbb{E}}{\widetilde{\bs{w}}}_i&={\mc{B}}\,{\mathbb{E}}{\widetilde{\bs{w}}}_{i-1}
-{\mc{A}}_2^{\T}\left(I_{NM}-{\mc{M}}{\mc{R}}'\right)\cdot{\mathbb{E}}{\bs{v}}_{i-1}^{(w)}
-{\mc{A}}_2^{\T}{\mc{M}}\cdot{\mathbb{E}}{\bs{z}}_i-{\mathbb{E}}{\bs{v}}_i^{(\psi)}
\end{align}
where
\begin{align}
\label{eqn:bigBdef}
{\mc{B}}&\defeq{\mc{A}}_2^{\T}\left(I_{NM}-{\mc{M}}{\mc{R}}'\right){\mc{A}}_1^{\T}\\
\label{eqn:bigRprimedef}
{\mc{R}}'&\defeq{\mathbb{E}}{\bs{\mc{R}}}_i'={\mathrm{diag}}\left\{R_1',\dots,R_N'\right\}\\
\label{eqn:Rkprimedef}
R_k'&\defeq{\mathbb{E}}{\bs{R}}_{k,i}'=\sum_{l\in{\mc{N}}_k}c_{lk}\left(R_{u,l}+R_{v,lk}^{(u)}\right)
\end{align}
From \eqref{eqn:vkiwdef}, \eqref{eqn:vkipsidef}, \eqref{eqn:viwdef}, and \eqref{eqn:vipsidef}, it can be verified that
\begin{align}
\label{eqn:meanviwvipsi}
{\mathbb{E}}{\bs{v}}_{i-1}^{(w)}={\mathbb{E}}{\bs{v}}_i^{(\psi)}=0
\end{align}
whereas, from \eqref{eqn:zkidef} and Assumption \ref{asm:all}, we get
\begin{align}
\label{eqn:meanzkidef}
{\mathbb{E}}{\bs{z}}_{k,i}&={\mathbb{E}}\left[\sum_{l\in{\mc{N}}_k}c_{lk}({\bs{u}}_{l,i}\!+\!{\bs{v}}_{lk,i}^{(u)})^*
({\bs{v}}_l(i)\!+\!{\bs{v}}_{lk}^{(d)}(i)\!-\!{\bs{v}}_{lk,i}^{(u)}w^o)\right]\nonumber\\
{}&=-\left(\sum_{l\in{\mc{N}}_k}c_{lk}R_{v,lk}^{(u)}\right)w^o
\end{align}
Let us define an $NM\times NM$ matrix ${\mc{R}}_{v,c}^{(u)}$ that collects all covariance matrices $\{R_{v,lk}^{(u)}\}$, $k,l=1,\dots,N$, weighted by the corresponding combination coefficients $\{c_{lk}\}$, such that its $(k,l)$th $M\times M$ submatrix is $c_{lk}R_{v,lk}^{(u)}$. Note that ${\mc{R}}_{v,c}^{(u)}$ itself is \emph{not} a covariance matrix because $c_{kk}R_{v,kk}^{(u)}=0$ for all $k$. Then, from \eqref{eqn:zidef} and \eqref{eqn:meanzkidef}, we arrive at
\begin{align}
\label{eqn:zdef}
z\defeq{\mathbb{E}}{\bs{z}}_{i}=-{\mc{R}}_{v,c}^{(u)}\left({\mathds{1}}_N\otimes w^o\right)
\end{align}
Therefore, using \eqref{eqn:meanviwvipsi} and \eqref{eqn:zdef}, expression \eqref{eqn:meanrecursion} becomes
\begin{equation}
\label{eqn:meanrecursion1}
\boxed{
{\mathbb{E}}{\wt{\bs{w}}}_i={\mc{B}}\cdot{\mathbb{E}}{\wt{\bs{w}}}_{i-1}-{\mc{A}}_2^{\T}{\mc{M}}z
}\end{equation}
with a driving term due to the presence of $z$. This driving term would disappear from \eqref{eqn:meanrecursion1} if there were no noise during the exchange of the regression data. To guarantee convergence of \eqref{eqn:meanrecursion1}, the coefficient matrix ${\mc{B}}$ must be stable, i.e., $\rho({\mc{B}})<1$. Since ${\mc{A}}_1^{\T}$ and ${\mc{A}}_2^{\T}$ are right-stochastic matrices, it can be shown that the matrix ${\mc{B}}$ is stable whenever $I_{NM}-{\mc{M}}{\mc{R}}'$ itself is stable (see Appendix \ref{app:meanconvergence}). This fact leads to an upper bound on the step-sizes $\{\mu_k\}$ to guarantee the convergence of ${\mathbb{E}}{\wt{\bs{w}}}_i$ to a steady-state value, namely, we must have
\begin{equation}
\boxed{
\label{eqn:meanconvergencecond2}
\mu_k<\frac{2}{\lambda_{\mathrm{max}}\left(R_k'\right)}
}\end{equation}
for $k=1,2,\dots,N$, where $\lambda_{\mathrm{max}}(\cdot)$ denotes the largest eigenvalue of its matrix argument. Note that the neighborhood covariance matrix $R_k'$ in \eqref{eqn:Rkprimedef} is related to the combination weights $\{c_{lk}\}$. If we further assume that $C$ is doubly-stochastic, i.e.,
\begin{align}
\label{eqn:doublyC}
C{\mathds{1}}_N={\mathds{1}}_N,\qquad C^{\T}{\mathds{1}}_N={\mathds{1}}_N
\end{align}
then, by Jensen's inequality \cite{Boyd04},
\begin{align}
\lambda_{\mathrm{max}}(R_k')&=\lambda_{\mathrm{max}}\left(\sum_{l\in{\mc{N}}_k}c_{lk}(R_{u,l}+R_{v,lk}^{(u)})\right)\nonumber\\
{}&\le\sum_{l\in{\mc{N}}_k}c_{lk}\lambda_{\mathrm{max}}\left(R_{u,l}+R_{v,lk}^{(u)}\right)\nonumber\\
{}&\le\max_{l\in\N_k}\lambda_{\mathrm{max}}\left(R_{u,l}+R_{v,lk}^{(u)}\right)
\end{align}
since (i) $\lambda_{\mathrm{max}}(\cdot)$ coincides with the induced $2$-norm for any positive semi-definite Hermitian matrix; (ii) matrix norms are convex functions of their arguments \cite{Horn85}; and (iii) by \eqref{eqn:doublyC}, $\{c_{lk}\}$ are convex combination coefficients. Thus, we obtain a sufficient condition for the convergence of \eqref{eqn:meanrecursion1} in lieu of \eqref{eqn:meanconvergencecond2}:
\begin{equation}
\label{eqn:meanconvergencecond3}
\boxed{
\mu_k<\frac{2}{\max_{l\in\N_k}\left[\lambda_{\mathrm{max}}\left(R_{u,l}+R_{v,lk}^{(u)}\right)\right]}
}\end{equation}
for $k=1,2,\dots,N$, where the upper bound for the step-size $\mu_k$ becomes independent of the combination weights $\{c_{lk}\}$. This bound can be determined solely from knowledge of the covariances of the regression data and the associated noise signals that are accessible to node $k$. It is worth noting that for traditional diffusion algorithms where information is perfectly exchanged, condition \eqref{eqn:meanconvergencecond3} reduces to
\begin{align}
\label{eqn:meanconvergencecond4}
\mu_{k}<\frac{2}{\max_{l\in\N_k}\left[\lambda_{\mathrm{max}}(R_{u,l})\right]}
\end{align}
for $k=1,2,\dots,N$. Comparing \eqref{eqn:meanconvergencecond3} with \eqref{eqn:meanconvergencecond4}, we see that the link noise ${\bs{v}}_{lk,i}^{(u)}$ over regression data reduces the dynamic range of the step-sizes for mean stability. Now, under \eqref{eqn:meanconvergencecond2}, and taking the limit of \eqref{eqn:meanrecursion1} as $i\rightarrow\infty$, we find that the mean error vector will converge to a steady-state value $g$:
\begin{align}
\label{eqn:gdef}
g\defeq\lim_{i\rightarrow\infty}{\mathbb{E}}{\widetilde{\bs{w}}}_i=-\left(I_{NM}-{\mc{B}}\right)^{-1}
{\mc{A}}_2^{\T}{\mc{M}}z
\end{align}

\section{Mean-Square Convergence Analysis}
It is well-known that studying the mean-square convergence of a single adaptive filter is a challenging task, since adaptive filters are nonlinear, time-variant, and stochastic systems. When a network of adaptive nodes is considered, the complexity of the analysis is compounded because the nodes now influence each other's behavior. In order to make the performance analysis more tractable, we rely on the energy conservation approach \cite{Sayed08,Alnaffouri03TSP}, which was used successfully in \cite{Lopes08TSP, Cattivelli10TSP} to study the mean-square performance of diffusion strategies under perfect information exchange conditions. That argument allows us to derive expressions for the mean-square-deviation (MSD) and the excess-mean-square-error (EMSE) of the network by analyzing how energy (measured in terms of error variances) flows through the nodes.

From recursion \eqref{eqn:noisyerrorrecursion1} and under Assumption \ref{asm:all}, we can obtain the following weighted variance relation for the global error vector ${\wt{\bs{w}}}_i$:
\begin{equation}
\boxed{
\begin{aligned}
\label{eqn:weightederrorrecursion}
{}\E\|{\wt{\bs{w}}}_i\|_{\Sigma}^2&=\E\|{\wt{\bs{w}}}_{i-1}\|_{\Sigma'}^2+\E\|{\mc{A}}_2^{\T}{\mc{M}}{\bs{z}}_i\|_{\Sigma}^2
-2\,\Re\{\E[{\bs{z}}_i^*{\mc{M}}{\mc{A}}_2\Sigma{\mc{A}}_2^{\T}(I_{NM}\!-\!{\mc{M}}{\bs{\mc{R}}}_i'){\mc{A}}_1^{\T}{\wt{\bs{w}}}_{i-1}]\}\\
{}&\quad+\E\|{\mc{A}}_2^{\T}(I_{NM}\!-\!{\mc{M}}{\bs{\mc{R}}}_i'){\bs{v}}_{i-1}^{(w)}\|_{\Sigma}^2+\E\|{\bs{v}}_i^{(\psi)}\|_{\Sigma}^2
\end{aligned}
}\end{equation}
where $\Sigma$ is an arbitrary $NM\times NM$ positive semi-definite Hermitian matrix that we are free to choose. Moreover, the notation $\|x\|_{\Sigma}^2$ stands for the quadratic term $x^*\Sigma x$. The weighting matrix $\Sigma'$ in \eqref{eqn:weightederrorrecursion} can be expressed as
\begin{align}
\label{eqn:Sigmaprimedef}
\Sigma'={\mc{B}}^*\Sigma{\mc{B}}+O({\mc{M}}^2)
\end{align}
where ${\mc{B}}$ is given by \eqref{eqn:bigBdef} and $O({\mc{M}}^2)$ denotes a term on the order of ${\mc{M}}^2$. Evaluating the term $O({\mc{M}}^2)$ requires knowledge of higher-order statistics of the regression data and link noises, which are not available under current assumptions. However, this term becomes negligible if we introduce a small step-size assumption.
\vspace{0.5\baselineskip}
\begin{assumption}[Small step-sizes]
\label{asm:smallstepsize}
The step-sizes are sufficiently small, i.e., $\mu_k\ll1$, such that terms depending on higher-order powers of the step-sizes can be ignored. \hfill \IEEEQED
\end{assumption}
\vspace{0.5\baselineskip}
Hence, in the sequel we use the approximation:
\begin{align}
\label{eqn:Sigmaprimedefsmall}
\Sigma'\approx{\mc{B}}^*\Sigma{\mc{B}}
\end{align}
Observe that on the right-hand side (RHS) of relation \eqref{eqn:weightederrorrecursion}, only the first and third terms relate to the error vector ${\wt{\bs{w}}}_{i-1}$. By Assumption \ref{asm:all}, the error vector ${\wt{\bs{w}}}_{i-1}$ is independent of ${\bs{z}}_i$ and ${\bs{\mc{R}}}_i'$. Thus, from \eqref{eqn:zdef}, the third term on RHS of \eqref{eqn:weightederrorrecursion} can be expressed as
\begin{align}
\label{eqn:thirdterm}
\mbox{Third term on RHS of \eqref{eqn:weightederrorrecursion}}&=-2\,\Re\{\E[{\bs{z}}_i^*{\mc{M}}{\mc{A}}_2\Sigma{\mc{A}}_2^{\T}(I_{NM}\!-\!{\mc{M}}{\bs{\mc{R}}}_i')
{\mc{A}}_1^{\T}]\cdot\E{\wt{\bs{w}}}_{i-1}\}\nonumber\\
{}&=-2\,\Re(z^*{\mc{M}}{\mc{A}}_2\Sigma{\mc{A}}_2^{\T}{\mc{A}}_1^{\T}\cdot\E{\wt{\bs{w}}}_{i-1})+O({\mc{M}}^2)
\end{align}
Since we already showed in the previous section that $\E{\wt{\bs{w}}}_{i}$ converges to a fixed bias $g$, quantity \eqref{eqn:thirdterm} will converge to a fixed value as well when $i\rightarrow\infty$. Moreover, under Assumption \ref{asm:all}, the second, fourth, and fifth terms on RHS of relation \eqref{eqn:weightederrorrecursion} are all fixed values. Therefore, the convergence of relation \eqref{eqn:weightederrorrecursion} depends on the behavior of the first term $\E\|{\wt{\bs{w}}}_{i-1}\|_{\Sigma'}^2$. Although the weighting matrix $\Sigma'$ of ${\wt{\bs{w}}}_{i-1}$ is different from the weighting matrix $\Sigma$ of ${\wt{\bs{w}}}_{i}$, it turns out that the entries of these two matrices are approximately related by a linear equation shown ahead in \eqref{eqn:linearrelation}. Introduce the vector notation \cite{Sayed08}:
\begin{align}
\label{eqn:vectornotation}
\sigma=\vec(\Sigma),\qquad\sigma'=\vec(\Sigma')
\end{align}
Then, by using the identity $\vec(ABC)=(C^{\T}\otimes A)\cdot\vec(B)$, it can be verified from \eqref{eqn:Sigmaprimedefsmall} that
\begin{align}
\label{eqn:linearrelation}
\sigma'\approx{\mc{F}}\cdot\sigma
\end{align}
where the $N^2M^2\times N^2M^2$ matrix ${\mc{F}}$ is given by
\begin{align}
\label{eqn:Fdefsmall}
{\mc{F}}&\defeq{\mc{B}}^{\T}\otimes{\mc{B}}^*
\end{align}
To guarantee mean-square convergence of the algorithm, the step-sizes should be sufficiently small and selected to ensure that the matrix ${\mc{F}}$ is stable \cite{Sayed08}, i.e., $\rho({\mc{F}})<1$, which is equivalent to the earlier condition $\rho({\mc{B}})<1$. Although more specific conditions for mean-square stability can be determined without Assumption \ref{asm:smallstepsize} \cite{Sayed08}, it is sufficient for our purposes here to conclude that the diffusion strategy \eqref{eqn:noisydiffusionpriordiffold}--\eqref{eqn:noisydiffusionpostdiffold} is stable in the mean and mean-square senses if the step-sizes $\{\mu_k\}$ satisfy \eqref{eqn:meanconvergencecond2} or \eqref{eqn:meanconvergencecond3} and are sufficiently small.

\section{Steady-State Performance Analysis}
The conclusion so far is that sufficiently small step-sizes ensure convergence of the diffusion strategy \eqref{eqn:noisydiffusionpriordiffold}--\eqref{eqn:noisydiffusionpostdiffold} in the mean and mean-square senses, even in the presence of exchange noises over the communication links. Let us now determine expressions for the error variances in steady-state. We start from the weighted variance relation \eqref{eqn:weightederrorrecursion}. In view of \eqref{eqn:thirdterm}, it shows that the error variance ${\mathbb{E}}\|{\wt{\bs{w}}}_i\|_{\Sigma}^2$ depends on the mean error ${\mathbb{E}}{\wt{\bs{w}}}_i$. We already determined the value of $\lim_{i\rightarrow\infty}{\mathbb{E}}{\wt{\bs{w}}}_i$ in \eqref{eqn:gdef}.

\subsection{Steady-State Variance Relation}
We continue to use the vector notation \eqref{eqn:vectornotation} and proceed to evaluate all the terms, except the first one, on RHS of \eqref{eqn:weightederrorrecursion} in the following. For the \emph{second} term, it can be expressed as
\begin{align}
\label{eqn:secondterm}
\mbox{Second term on RHS of \eqref{eqn:weightederrorrecursion}}&=\Tr({\mc{A}}_2^{\T}{\mc{M}}{\mc{R}}_z{\mc{M}}{\mc{A}}_2\Sigma)\nonumber\\
&=\left[\vec({\mc{A}}_2^{\T}{\mc{M}}{\mc{R}}_z{\mc{M}}{\mc{A}}_2)\right]^*\sigma
\end{align}
where we used the identity $\Tr(W\Sigma)=[\vec(W)]^*\sigma$ for any Hermitian matrix $W$, and ${\mc{R}}_{z}$ denotes the autocorrelation matrix of ${\bs{z}}_{i}$. It is shown in Appendix \ref{app:Rz} that ${\mc{R}}_{z}$ is given by
\begin{align}
\label{eqn:bigRzdef}
{\mc{R}}_z\defeq\E\,{\bs{z}}_i{\bs{z}}_i^*\approx{\mc{C}}^{\T}{\mc{S}}{\mc{C}}+{\mc{T}}+zz^*
\end{align}
where ${\mc{C}}$ is defined in \eqref{eqn:bigA1CA2def}, $z$ is in \eqref{eqn:zdef},  and $\{{\mc{S}},{\mc{T}}\}$ are two $NM\times NM$ positive semi-definite block diagonal matrices:
\begin{align}
\label{eqn:bigSdef}
{\mc{S}}&\defeq{\mathrm{diag}}\left\{\sigma_{v,1}^{2}R_{u,1},\dots,\sigma_{v,N}^{2}R_{u,N}\right\}\\
\label{eqn:bigTdef}
{\mc{T}}&\defeq{\mathrm{diag}}\left\{T_{1},\dots,T_{N}\right\}\\
\label{eqn:Tkdef}
T_{k}&\defeq\sum_{l\in{\mc{N}}_k}\!c_{lk}^2\left[(\sigma_{v,l}^{2}\!+\!\sigma_{v,lk}^{2})R_{v,lk}^{(u)}
\!+\!(\sigma_{v,lk}^{2}\!+\!w^{o*}R_{v,lk}^{(u)}w^o)R_{u,l}\right]
\end{align}
From expression \eqref{eqn:thirdterm} and Assumption \ref{asm:smallstepsize}, the \emph{third} term on RHS of \eqref{eqn:weightederrorrecursion} is given by
\begin{align}
\mbox{Third term on RHS of \eqref{eqn:weightederrorrecursion}}&
\approx-z^*{\mc{M}}{\mc{A}}_2\Sigma{\mc{A}}_2^{\T}{\mc{A}}_1^{\T}(\E{\wt{\bs{w}}}_{i-1})
-(\E{\wt{\bs{w}}}_{i-1})^*{\mc{A}}_1{\mc{A}}_2\Sigma{\mc{A}}_2^{\T}{\mc{M}}z\nonumber\\
{}&=-\Tr\{[{\mc{A}}_2^{\T}{\mc{A}}_1^{\T}(\E{\wt{\bs{w}}}_{i-1})z^*\!{\mc{M}}{\mc{A}}_2
\!+\!{\mc{A}}_2^{\T}{\mc{M}}z(\E{\wt{\bs{w}}}_{i-1})^*\!{\mc{A}}_1{\mc{A}}_2]\Sigma\}\nonumber\\
{}&=-\!\left[\vec\left({\mc{A}}_2^{\T}{\mc{A}}_1^{\T}(\E{\wt{\bs{w}}}_{i-1})z^*{\mc{M}}{\mc{A}}_2
+{\mc{A}}_2^{\T}{\mc{M}}z(\E{\wt{\bs{w}}}_{i-1})^*{\mc{A}}_1{\mc{A}}_2\right)\right]^*\sigma
\end{align}
Likewise, the \emph{fourth} term on RHS of \eqref{eqn:weightederrorrecursion} is approximated by
\begin{align}
\mbox{Fourth term on RHS of \eqref{eqn:weightederrorrecursion}}
&=\left[\vec\!\left(\E[{\mc{A}}_2^{\T}(I_{NM}\!-\!{\mc{M}}{\bs{\mc{R}}}_i'){\mc{R}}_v^{(w)}
(I_{NM}\!-\!{\mc{M}}{\bs{\mc{R}}}_i'){\mc{A}}_2]\right)\right]^*\sigma\nonumber\\
{}&\approx\left[\vec({\mc{A}}_2^{\T}{\mc{R}}_v^{(w)}{\mc{A}}_2)\right]^*\sigma
\end{align}
where we are now ignoring terms on the order of ${\mc{M}}$ and ${\mc{M}}^2$. The \emph{fifth} term on RHS of \eqref{eqn:weightederrorrecursion} is given by
\begin{align}
\mbox{Fifth term on RHS of \eqref{eqn:weightederrorrecursion}}=\left[\vec({\mc{R}}_v^{(\psi)})\right]^*\sigma
\end{align}
Let us introduce
\begin{align}
\label{eqn:bigRvdef}
{\mc{R}}_v&\defeq{\mc{A}}_2^{\T}{\mc{R}}_v^{(w)}{\mc{A}}_2+{\mc{R}}_v^{(\psi)}
+{\mc{A}}_2^{\T}{\mc{M}}({\mc{T}}+zz^*){\mc{M}}{\mc{A}}_2\\
\label{eqn:bigYdef}
{\mc{Y}}&\defeq-{\mc{A}}_2^{\T}{\mc{A}}_1^{\T}gz^*{\mc{M}}{\mc{A}}_2\nonumber\\
{}&={\mc{A}}_2^{\T}{\mc{A}}_1^{\T}\left(I_{NM}-{\mc{B}}\right)^{-1}{\mc{A}}_2^{\T}{\mc{M}}zz^*{\mc{M}}{\mc{A}}_2
\end{align}
At steady-state, as $i\rightarrow\infty$, by \eqref{eqn:gdef} and \eqref{eqn:secondterm}--\eqref{eqn:bigYdef},  the weighted variance relation \eqref{eqn:weightederrorrecursion} becomes
\begin{align}
\label{eqn:weightederrorrecursion1}
\lim_{i\rightarrow\infty}&{\mathbb{E}}\|{\wt{\bs{w}}}_{i}\|_{\sigma}^2\approx
\lim_{i\rightarrow\infty}{\mathbb{E}}\|{\wt{\bs{w}}}_{i-1}\|_{{\mc{F}}\sigma}^2
+\left[\vec({\mc{A}}_2^{\T}{\mc{M}}{\mc{C}}^{\T}{\mc{S}}{\mc{C}}{\mc{M}}{\mc{A}}_2
+{\mc{R}}_v+{\mc{Y}}+{\mc{Y}}^*)\right]^*\sigma
\end{align}
where we are using the compact notation $\|x\|^2_{\sigma}$ to refer to $\|x\|^2_{\Sigma}$ --- doing so allows us to represent $\Sigma'$ by the more compact relation ${\mc{F}}\sigma$ on RHS of \eqref{eqn:weightederrorrecursion1}; we shall be using the weighting matrix $\Sigma$ and its vector representation $\sigma$ interchangeably for ease of notation (likewise, for $\Sigma'$ and $\sigma'$). The steady-state weighted variance relation \eqref{eqn:weightederrorrecursion1} can be rewritten as
\begin{align}
\label{eqn:steadystateweightednew}
\lim_{i\rightarrow\infty}&{\mathbb{E}}\|{\wt{\bs{w}}}_{i}\|_{(I_{N^2M^2}-{\mc{F}})\sigma}^2
\approx\left[\vec({\mc{A}}_2^{\T}{\mc{M}}{\mc{C}}^{\T}{\mc{S}}{\mc{C}}{\mc{M}}{\mc{A}}_2
+{\mc{R}}_v+{\mc{Y}}+{\mc{Y}}^*)\right]^*\sigma
\end{align}
where the term ${\mc{A}}_2^{\T}{\mc{M}}{\mc{C}}^{\T}{\mc{S}}{\mc{C}}{\mc{M}}{\mc{A}}_2$ is contributed by the model noise $\{{\bs{v}}_k(i)\}$ while the remaining terms $\{{\mc{R}}_v,{\mc{Y}}\}$ are contributed by the link noises $\{{\bm{v}}_{lk,i-1}^{(w)},{\bm{v}}_{lk}^{(d)}(i),{\bm{v}}_{lk,i}^{(u)},{\bm{v}}_{lk,i}^{(\psi)}\}$. Recall that we are free to choose $\Sigma$ and, hence, $\sigma$. Let $(I_{N^2M^2}-{\mc{F}})\sigma={\mathrm{vec}}(\Omega)$, where $\Omega$ is another arbitrary positive semi-definite Hermitian matrix. Then, we arrive at the following theorem.
\vspace{0.5\baselineskip}
\begin{theorem}[Steady-state weighted variance relation]
\label{lemma:steadystatevariancerelation}
Under Assumptions \ref{asm:all} and \ref{asm:smallstepsize}, for any positive semi-definite Hermitian matrix $\Omega$, the steady-state weighted error variance relation of the diffusion strategy \eqref{eqn:noisydiffusionpriordiffold}--\eqref{eqn:noisydiffusionpostdiffold} is approximately given by
\begin{equation}
\label{eqn:steadystatevariance}
\boxed{
\begin{aligned}
\lim_{i\rightarrow\infty}{\mathbb{E}}\|{\wt{\bs{w}}}_{i}\|_{\Omega}^2\approx
\left[\vec({\mc{A}}_2^{\T}{\mc{M}}{\mc{C}}^{\T}{\mc{S}}{\mc{C}}{\mc{M}}{\mc{A}}_2
\!+\!{\mc{R}}_v\!+\!{\mc{Y}}\!+\!{\mc{Y}}^*)\right]^*(I_{N^2M^2}-{\mc{F}})^{-1}\vec(\Omega)
\end{aligned}
}\end{equation}
where ${\mc{S}}$ is given in \eqref{eqn:bigSdef}, ${\mc{R}}_v$ in \eqref{eqn:bigRvdef}, ${\mc{Y}}$ in \eqref{eqn:bigYdef}, and ${\mc{F}}$ in \eqref{eqn:Fdefsmall}. \hfill \IEEEQED
\end{theorem}

\subsection{Network MSD and EMSE}
Each subvector of ${\widetilde{\bs{w}}}_{i}$ corresponds to the estimation error at a particular node, say, ${\widetilde{\bs{w}}}_{k,i}$ for node $k$. The network MSD is defined as \cite{Sayed08}:
\begin{align}
{\textrm{MSD}}\defeq\lim_{i\rightarrow\infty}\frac{1}{N}\sum_{k=1}^{N}{\mathbb{E}}\|{\widetilde{\bs{w}}}_{k,i}\|^2
\end{align}
Since we are free to choose $\Omega$, we select it as $\Omega=I_{NM}/N$. Then, expression \eqref{eqn:steadystatevariance} gives
\begin{equation}
\boxed{
\label{eqn:noisyMSD}
\begin{aligned}
{\textrm{MSD}}&\!\approx\!\frac{1}{N}\!\left[\vec({\mc{A}}_2^{\T}{\mc{M}}{\mc{C}}^{\T}{\mc{S}}{\mc{C}}{\mc{M}}{\mc{A}}_2
\!+\!{\mc{R}}_v\!+\!{\mc{Y}}\!+\!{\mc{Y}}^*)\right]^*(I_{N^2M^2}-{\mc{F}})^{-1}\vec(I_{NM})
\end{aligned}
}\end{equation}
Similarly, if we instead select $\Omega={\mc{R}}_u/N$, where
\begin{align}
\label{eqn:bigRudef}
{\mc{R}}_u\defeq{\mathrm{diag}}\left\{R_{u,1},\dots,R_{u,N}\right\}
\end{align}
then expression \eqref{eqn:steadystatevariance} would allow us to evaluate the network EMSE as:
\begin{equation}
\boxed{
\label{eqn:noisyEMSE}
\begin{aligned}
{\textrm{EMSE}}&\!\approx\!\frac{1}{N}\!\left[\vec({\mc{A}}_2^{\T}{\mc{M}}{\mc{C}}^{\T}{\mc{S}}{\mc{C}}{\mc{M}}{\mc{A}}_2
\!+\!{\mc{R}}_v\!+\!{\mc{Y}}\!+\!{\mc{Y}}^*)\right]^*(I_{N^2M^2}-{\mc{F}})^{-1}\vec({\mc{R}}_u)
\end{aligned}
}\end{equation}
where, under Assumption \ref{asm:all}, the network EMSE is given by
\begin{align}
\label{eqn:EMSEdef}
{\textrm{EMSE}}&\defeq\lim_{i\rightarrow\infty}\frac{1}{N}\sum_{k=1}^{N}{\mathbb{E}}|{\bs{u}}_{k,i}{\wt{\bs{w}}}_{k,i-1}|^2\nonumber\\
{}&=\lim_{i\rightarrow\infty}\frac{1}{N}\sum_{k=1}^{N}{\mathbb{E}}\|{\wt{\bs{w}}}_{k,i}\|_{{\mc{R}}_u}^2
\end{align}

\subsection{Simplifications when Regression Data are not Shared}
We showed in the earlier sections that the link noise over regression data biases the weight estimators. In this section we examine how the results simplify when there is no sharing of regression data among the nodes.
\vspace{0.5\baselineskip}
\begin{assumption}[No sharing of regression data]
\label{asm:nodatasharing} Nodes do not share regression data within neighborhoods, i.e., assume $C=I_N$. \hfill \IEEEQED
\end{assumption}
\vspace{0.5\baselineskip}
By Assumptions \ref{asm:smallstepsize} and \ref{asm:nodatasharing}, matrices $\{{\mc{B}},{\mc{R}}_v,{\mc{Y}}\}$ in \eqref{eqn:bigBdef}, \eqref{eqn:bigRvdef}, and \eqref{eqn:bigYdef} become
\begin{align}
\label{eqn:bigBsimple}
{\mc{B}}&={\mc{A}}_2^{\T}(I_{NM}-{\mc{M}}{\mc{R}}_u){\mc{A}}_1^{\T}\\
\label{eqn:bigRvsimple}
{\mc{R}}_v&={\mc{A}}_2^{\T}{\mc{R}}_v^{(w)}{\mc{A}}_2+{\mc{R}}_v^{(\psi)}\\
\label{eqn:bigYsimple}
{\mc{Y}}&=0
\end{align}
where ${\mc{R}}_u$ is given in \eqref{eqn:bigRudef}. Then, the network MSD and EMSE expressions \eqref{eqn:noisyMSD} and \eqref{eqn:noisyEMSE} simplify to:
\begin{equation}
\boxed{
\label{eqn:noisyMSD1}
\begin{aligned}
{\textrm{MSD}}&\approx\frac{1}{N}\left[\vec({\mc{A}}_2^{\T}{\mc{M}}{\mc{S}}{\mc{M}}{\mc{A}}_2+{\mc{R}}_v)\right]^*
(I_{N^2M^2}-{\mc{F}})^{-1}\vec(I_{NM})
\end{aligned}
}\end{equation}
and
\begin{equation}
\boxed{
\label{eqn:noisyEMSE1}
\begin{aligned}
{\textrm{EMSE}}&\approx\frac{1}{N}\left[\vec({\mc{A}}_2^{\T}{\mc{M}}{\mc{S}}{\mc{M}}{\mc{A}}_2+{\mc{R}}_v)\right]^*
(I_{N^2M^2}-{\mc{F}})^{-1}\vec({\mc{R}}_u)
\end{aligned}
}\end{equation}

\subsection{Dependence of Performance on Combination Weights and Link Noise}
Recalling that ${\mc{R}}_v$ and ${\mc{F}}$ are related to the combination matrices $\{{\mc{A}}_1,{\mc{A}}_2\}$, or, equivalently, $\{A_1,A_2\}$, results \eqref{eqn:noisyMSD1} and \eqref{eqn:noisyEMSE1} express the network MSD and EMSE in terms of $\{A_1,A_2\}$. However, it is generally difficult to use these expressions to optimize over $\{A_1,A_2\}$ to reduce the impact of link noise. Instead, by substituting \eqref{eqn:Fdefsmall} into \eqref{eqn:noisyMSD1} and using the fact that ${\mc{F}}$ is stable, we can arrive at another useful expression for the network MSD:
\begin{align}
{\textrm{MSD}}&\approx\frac{1}{N}\!\left[\vec({\mc{A}}_2^{\T}{\mc{M}}{\mc{S}}{\mc{M}}{\mc{A}}_2\!+\!{\mc{R}}_v)\right]^*
\sum_{j=0}^{\infty}{\mc{F}}^j\vec(I_{NM})\nonumber\\
{}&=\frac{1}{N}\!\left[\vec({\mc{A}}_2^{\T}{\mc{M}}{\mc{S}}{\mc{M}}{\mc{A}}_2\!+\!{\mc{R}}_v)\right]^*
\sum_{j=0}^{\infty}({\mc{B}}^{\T}\otimes{\mc{B}}^*)^j\vec(I_{NM})\nonumber\\
{}&=\frac{1}{N}\!\left[\vec({\mc{A}}_2^{\T}{\mc{M}}{\mc{S}}{\mc{M}}{\mc{A}}_2\!+\!{\mc{R}}_v)\right]^*
\sum_{j=0}^{\infty}\vec({\mc{B}}^{*j}{\mc{B}}^j)
\end{align}
That is,
\begin{equation}
\label{eqn:noisyMSD2}
\boxed{
{\textrm{MSD}}\approx\frac{1}{N}\sum_{j=0}^{\infty}{\mathrm{Tr}}
\left[{\mc{B}}^j({\mc{A}}_2^{\T}{\mc{M}}{\mc{S}}{\mc{M}}{\mc{A}}_2\!+\!{\mc{R}}_v){\mc{B}}^{*j}\right]
}\end{equation}
where ${\mc{B}}$ is given in \eqref{eqn:bigBsimple}. Similarly, the network EMSE can be expressed as
\begin{equation}
\label{eqn:noisyEMSE2}
\boxed{
{\textrm{EMSE}}\approx\frac{1}{N}\sum_{j=0}^{\infty}{\mathrm{Tr}}
\left[{\mc{B}}^j({\mc{A}}_2^{\T}{\mc{M}}{\mc{S}}{\mc{M}}{\mc{A}}_2\!+\!{\mc{R}}_v){\mc{B}}^{*j}{\mc{R}}_u\right]
}\end{equation}
Expressions \eqref{eqn:noisyMSD2} and \eqref{eqn:noisyEMSE2} reveal in an interesting way how the noise sources originating from any particular node end up influencing the overall network performance. Let us denote
\begin{align}
\label{eqn:rndbigBidef}
{\bs{\mc{B}}}_i&\defeq{\mc{A}}_2^{\T}(I_{NM}-{\mc{M}}{\bs{\mc{R}}}_i'){\mc{A}}_1^{\T}\\
\label{eqn:thetaidef}
{\bs{\theta}}_i&\defeq{\mc{A}}_2^{\T}(I_{NM}-{\mc{M}}{\bs{\mc{R}}}_i'){\bs{v}}_{i-1}^{(w)}
+{\mc{A}}_2^{\T}{\mc{M}}{\bs{z}}_i+{\bs{v}}_i^{(\psi)}
\end{align}
The error recursion \eqref{eqn:noisyerrorrecursion1} can be rewritten as
\begin{align}
\label{eqn:noisyerrorrecursion2}
{\wt{\bs{w}}}_i&={\bs{\mc{B}}}_i{\wt{\bs{w}}}_{i-1}-{\bs{\theta}}_i\nonumber\\
{}&={\bm{\Phi}}_{0,i}{\wt{\bs{w}}}_{-1}-\sum_{m=0}^{i}{\bm{\Phi}}_{m+1,i}{\bs{\theta}}_m
\end{align}
where
\begin{align}
{\bm{\Phi}}_{m,i}\defeq\begin{cases}
{\bs{\mc{B}}}_{i}{\bs{\mc{B}}}_{i-1}\dots{\bs{\mc{B}}}_{m}, & i\ge m \\
I_{NM}, & i<m
\end{cases}
\end{align}
Then,
\begin{align}
\label{eqn:variancerelationnew}
\E\,\|{\wt{\bs{w}}}_i\|^2=\E\,\|{\bm{\Phi}}_{0,i}{\wt{\bs{w}}}_{-1}\|^2+\E\left\|\sum_{m=0}^{i}{\bm{\Phi}}_{m+1,i}{\bs{\theta}}_m\right\|^2
\end{align}
Under Assumption \ref{asm:nodatasharing}, $\{{\bs{\mc{B}}}_i,{\bs{\theta}}_i\}$ in \eqref{eqn:rndbigBidef} and \eqref{eqn:thetaidef} can be simplified as
\begin{align}
\label{eqn:rndbigBisimple}
{\bs{\mc{B}}}_i&={\mc{A}}_2^{\T}(I_{NM}-{\mc{M}}{\bs{\mc{R}}}_i){\mc{A}}_1^{\T}\\
\label{eqn:thetaisimple}
{\bs{\theta}}_i&={\mc{A}}_2^{\T}(I_{NM}-{\mc{M}}{\bs{\mc{R}}}_i){\bs{v}}_{i-1}^{(w)}
+{\mc{A}}_2^{\T}{\mc{M}}{\bs{s}}_i+{\bs{v}}_i^{(\psi)}
\end{align}
where $\{{\bs{\mc{R}}}_i,{\bs{s}}_i\}$ are given in \eqref{eqn:rndbigRidef} and \eqref{eqn:sidef}. By Assumption \ref{asm:all}, $\{{\bs{\mc{B}}}_i,{\bs{\theta}}_i\}$ are temporally independent for different $i$ and
\begin{align}
\E\,{\bs{\mc{B}}}_i={\mc{B}},\qquad \E\,{\bs{\theta}}_i=0
\end{align}
where ${\mc{B}}$ is given by \eqref{eqn:bigBsimple}. As $i\rightarrow\infty$, the first term on RHS of \eqref{eqn:variancerelationnew} becomes
\begin{align}
\label{eqn:initialterm}
\mbox{First term on RHS of \eqref{eqn:variancerelationnew}}&
=\lim_{i\rightarrow\infty}\Tr\left\{\E\left[{\bm{\Phi}}_{0,i}(\E{\wt{\bs{w}}}_{-1}{\wt{\bs{w}}}_{-1}^*){\bm{\Phi}}_{0,i}^*\right]\right\}\nonumber\\
{}&\stackrel{(a)}{\approx}\lim_{i\rightarrow\infty}\Tr\left[\left(\E{\bm{\Phi}}_{0,i}\right)
(\E{\wt{\bs{w}}}_{-1}{\wt{\bs{w}}}_{-1}^*)\left(\E{\bm{\Phi}}_{0,i}\right)^*\right]\nonumber\\
{}&=\lim_{i\rightarrow\infty}\Tr\left[{\mc{B}}^{i+1}(\E\,{\wt{\bs{w}}}_{-1}{\wt{\bs{w}}}_{-1}^*){\mc{B}}^{(i+1)*}\right]\nonumber\\
{}&\stackrel{(b)}{=}0
\end{align}
where (a) is obtained by approximating the expectation of the product by the product of expectations and (b) is due to the stability of ${\mc{B}}$. Therefore, the steady-state value of \eqref{eqn:variancerelationnew} gives
\begin{align}
\label{eqn:steadystatevariancenew1}
\lim_{i\rightarrow\infty}\E\,\|{\wt{\bs{w}}}_i\|^2&=
\lim_{i\rightarrow\infty}\sum_{m=0}^{i}\E\,\left\|{\bm{\Phi}}_{m+1,i}{\bs{\theta}}_m\right\|^2\nonumber\\
{}&\approx\lim_{i\rightarrow\infty}\sum_{m=0}^{i}\Tr\left[\left(\E{\bm{\Phi}}_{m+1,i}\right)
(\E{\bs{\theta}}_m{\bs{\theta}}_m^*)\left(\E{\bm{\Phi}}_{m+1,i}\right)^*\right]\nonumber\\
{}&\stackrel{(a)}{\approx}\lim_{i\rightarrow\infty}\sum_{m=0}^{i}\Tr\left[{\mc{B}}^{i-m}({\mc{A}}_2^{\T}{\mc{M}}{\mc{S}}
{\mc{M}}{\mc{A}}_2+{\mc{R}}_v){\mc{B}}^{(i-m)*}\right]\nonumber\\
{}&\stackrel{(b)}{=}\lim_{i\rightarrow\infty}\sum_{j=0}^{i}\Tr\left[{\mc{B}}^{j}({\mc{A}}_2^{\T}{\mc{M}}{\mc{S}}
{\mc{M}}{\mc{A}}_2+{\mc{R}}_v){\mc{B}}^{j*}\right]\nonumber\\
{}&=\sum_{j=0}^{\infty}\Tr\left[{\mc{B}}^{j}({\mc{A}}_2^{\T}{\mc{M}}{\mc{S}}{\mc{M}}{\mc{A}}_2+{\mc{R}}_v){\mc{B}}^{*j}\right]
\end{align}
where, by \eqref{eqn:bigRvsimple} and \eqref{eqn:thetaisimple}, (a) is due to
\begin{align}
\E{\bs{\theta}}_m{\bs{\theta}}_m^*&\approx{\mc{A}}_2^{\T}(I_{NM}-{\mc{M}}{\mc{R}}_u){\mc{R}}_{v}^{(w)}(I_{NM}-{\mc{M}}{\mc{R}}_u){\mc{A}}_2
+{\mc{A}}_2^{\T}{\mc{M}}{\mc{S}}{\mc{M}}{\mc{A}}_2+{\mc{R}}_{v}^{(\psi)}\nonumber\\
{}&\approx{\mc{A}}_2^{\T}{\mc{M}}{\mc{S}}{\mc{M}}{\mc{A}}_2+{\mc{A}}_2^{\T}{\mc{R}}_{v}^{(w)}{\mc{A}}_2+{\mc{R}}_{v}^{(\psi)}\nonumber\\
{}&={\mc{A}}_2^{\T}{\mc{M}}{\mc{S}}{\mc{M}}{\mc{A}}_2+{\mc{R}}_{v}
\end{align}
and (b) is simply a change of variable: $j=i-m$. Since the $j$th term of the summation in \eqref{eqn:noisyMSD2} or \eqref{eqn:steadystatevariancenew1} is contributed by the term $\E{\bs{\theta}}_{i-j}{\bs{\theta}}_{i-j}^*$, which consists of all the noise sources at time $i-j$, expression \eqref{eqn:noisyMSD2} shows how various sources of noises are involved and how they contribute to the network MSD.

\section{Optimizing the Combination Matrices}
Before we optimize the combination matrices $\{A_1,A_2\}$, we first specialize the MSD expression  \eqref{eqn:noisyMSD2} and the EMSE expression \eqref{eqn:noisyEMSE2} for the ATC and CTA algorithms. For the ATC algorithm, we set $A_1=I_N$ and $A_2=A$, and for the CTA algorithm, we set $A_1=A$ and $A_2=I_N$. Let us denote
\begin{align}
{\mc{A}}&\defeq A\otimes I_M\\
{\mc{B}}_{\textrm{atc}}&\defeq{\mc{A}}^{\T}(I_{NM}-{\mc{M}}{\mc{R}}_u)\\
{\mc{B}}_{\textrm{cta}}&\defeq(I_{NM}-{\mc{M}}{\mc{R}}_u){\mc{A}}^{\T}
\end{align}
Then, we get
\begin{align}
\label{eqn:noisyatcMSD}
\!\!{\textrm{MSD}}_{\textrm{atc}}&\!\approx\!\frac{1}{N}\!\sum_{j=0}^{\infty}{\mathrm{Tr}}\!\left[{\mc{B}}_{\textrm{atc}}^j
({\mc{A}}^{\T}{\mc{M}}{\mc{S}}{\mc{M}}{\mc{A}}\!+\!{\mc{R}}_v^{(\psi)})
{\mc{B}}_{\textrm{atc}}^{*j}\right]\!\!\\
\label{eqn:noisyatcEMSE}
\!\!{\textrm{EMSE}}_{\textrm{atc}}&\!\approx\!\frac{1}{N}\!\sum_{j=0}^{\infty}{\mathrm{Tr}}\!\left[{\mc{B}}_{\textrm{atc}}^j
({\mc{A}}^{\T}{\mc{M}}{\mc{S}}{\mc{M}}{\mc{A}}\!+\!{\mc{R}}_v^{(\psi)})
{\mc{B}}_{\textrm{atc}}^{*j}{\mc{R}}_u\right]\!\!
\end{align}
and
\begin{align}
\label{eqn:noisyctaMSD}
{\textrm{MSD}}_{\textrm{cta}}&\approx\frac{1}{N}\sum_{j=0}^{\infty}{\mathrm{Tr}}\left[{\mc{B}}_{\textrm{cta}}^j
({\mc{M}}{\mc{S}}{\mc{M}}+{\mc{R}}_v^{(w)}){\mc{B}}_{\textrm{cta}}^{*j}\right]\\
\label{eqn:noisyctaEMSE}
{\textrm{EMSE}}_{\textrm{cta}}&\approx\frac{1}{N}\sum_{j=0}^{\infty}{\mathrm{Tr}}\left[{\mc{B}}_{\textrm{cta}}^j
({\mc{M}}{\mc{S}}{\mc{M}}+{\mc{R}}_v^{(w)}){\mc{B}}_{\textrm{cta}}^{*j}{\mc{R}}_u\right]
\end{align}

\subsection{An Upper Bound for MSD}
Minimizing the MSD expression \eqref{eqn:noisyatcMSD} or the EMSE expression \eqref{eqn:noisyatcEMSE} for the ATC algorithm over left-stochastic matrices $A$ is generally nontrivial. We pursue an approximate solution that relies on optimizing an upper bound and performs well in practice. Let us use $\|X\|_*$ to denote the nuclear norm (also known as the trace norm, or the Ky Fan $n$-norm) of matrix $X$ \cite{Laub05}, which is defined as the sum of the singular values of $X$. Therefore, $\|X\|_*=\|X^*\|_*$ for any $X$ and $\|X\|_*={\mathrm{Tr}}(X)$ when $X$ is Hermitian and positive semi-definite. Let us also denote $\|X\|_{b,\infty}$ as the block maximum norm of matrix $X$ (see Appendix \ref{app:meanconvergence}). Then,
\begin{align}
\label{eqn:uptrace}
{\mathrm{Tr}}\left[{\mc{B}}_{\textrm{atc}}^j({\mc{A}}^{\T}{\mc{M}}{\mc{S}}{\mc{M}}
{\mc{A}}+{\mc{R}}_v^{(\psi)}){\mc{B}}_{\textrm{atc}}^{*j}\right]
&=\big\|{\mc{B}}_{\textrm{atc}}^j({\mc{A}}^{\T}{\mc{M}}{\mc{S}}{\mc{M}}{\mc{A}}
+{\mc{R}}_v^{(\psi)}){\mc{B}}_{\textrm{atc}}^{*j}\big\|_*\nonumber\\
{}&\le\|{\mc{B}}_{\textrm{atc}}^j\|_*\cdot\|{\mc{A}}^{\T}{\mc{M}}{\mc{S}}
{\mc{M}}{\mc{A}}+{\mc{R}}_v^{(\psi)}\|_*\cdot\|{\mc{B}}_{\textrm{atc}}^{*j}\|_*\nonumber\\
{}&\le c^2\cdot\|{\mc{B}}_{\textrm{atc}}^j\|_{b,\infty}^2\cdot
{\mathrm{Tr}}({\mc{A}}^{\T}{\mc{M}}{\mc{S}}{\mc{M}}{\mc{A}}+{\mc{R}}_v^{(\psi)})\nonumber\\
{}&\le c^2\cdot\|{\mc{B}}_{\textrm{atc}}\|_{b,\infty}^{2j}\cdot
{\mathrm{Tr}}({\mc{A}}^{\T}{\mc{M}}{\mc{S}}{\mc{M}}{\mc{A}}+{\mc{R}}_v^{(\psi)})\nonumber\\
{}&\le c^2\cdot(\|{\mc{A}}\|_{b,\infty}\cdot\|I_{NM}-{\mc{M}}{\mc{R}}_u\|_{b,\infty})^{2j}
{\mathrm{Tr}}({\mc{A}}^{\T}{\mc{M}}{\mc{S}}{\mc{M}}{\mc{A}}+{\mc{R}}_v^{(\psi)})\nonumber\\
{}&=c^2\cdot\rho(I_{NM}-{\mc{M}}{\mc{R}}_u)^{2j}\cdot{\mathrm{Tr}}
({\mc{A}}^{\T}{\mc{M}}{\mc{S}}{\mc{M}}{\mc{A}}+{\mc{R}}_v^{(\psi)})
\end{align}
where $c$ is some positive scalar such that $\|X\|_*\le c\,\|X\|_{b,\infty}$ because $\|X\|_*$ and $\|X\|_{b,\infty}$ are submultiplicative norms and all such norms are equivalent \cite{Horn85}. In the last step of \eqref{eqn:uptrace} we used Lemmas \ref{lemma:rghtstocmat} and \ref{lemma:blockdiagonal} from Appendix \ref{app:meanconvergence}. Thus, we can upper bound the network MSD \eqref{eqn:noisyatcMSD} by
\begin{align}
\label{eqn:noisyatcMSDub}
{\textrm{MSD}}_{\textrm{atc}}&\le\frac{1}{N}\sum_{j=0}^{\infty}c^2\cdot\rho(I_{NM}-{\mc{M}}{\mc{R}}_u)^{2j}
\Tr({\mc{A}}^{\T}{\mc{M}}{\mc{S}}{\mc{M}}{\mc{A}}+{\mc{R}}_v^{(\psi)})\nonumber\\
{}&=\frac{c^2}{N}\cdot\frac{{\mathrm{Tr}}({\mc{A}}^{\T}{\mc{M}}{\mc{S}}{\mc{M}}{\mc{A}}+{\mc{R}}_v^{(\psi)})}
{1-[\rho(I_{NM}-{\mc{M}}{\mc{R}}_u)]^2}
\end{align}
where the combination matrix ${\mc{A}}$ appears only in the numerator.

\subsection{Minimizing the Upper Bound}
The result \eqref{eqn:noisyatcMSDub} motivates us to consider instead the problem of minimizing the upper bound, namely,
\begin{align}
\label{eqn:minproblem}
\begin{aligned}
\minimize_A \quad & \qquad \Tr({\mc{A}}^{\T}{\mc{M}}{\mc{S}}{\mc{M}}{\mc{A}}+{\mc{R}}_v^{(\psi)})\\
\st \quad & A^{\T}{\mathds{1}}={\mathds{1}},\;\;a_{lk}\ge0,\;\;a_{lk}=0\;\mbox{if}\;l\notin{\mc{N}}_k  \\
\end{aligned}
\end{align}
Using \eqref{eqn:bigRvpsidef} and \eqref{eqn:bigSdef}, the cost function in \eqref{eqn:minproblem} can be expressed as
\begin{align}
\Tr({\mc{A}}^{\T}{\mc{M}}{\mc{S}}{\mc{M}}{\mc{A}}+{\mc{R}}_v^{(\psi)})
=\sum_{k=1}^{N}\sum_{l\in{\mc{N}}_k}a_{lk}^2\left[\mu_l^2\sigma_{v,l}^2\Tr(R_{u,l})+\Tr(R_{v,lk}^{(\psi)})\right]
\end{align}
Problem \eqref{eqn:minproblem} can therefore be decoupled into $N$ separate optimization problems of the form:
\begin{equation}
\boxed{
\label{eqn:minsubproblem}
\begin{aligned}
\minimize_{\{a_{lk},\;l\in{\mc{N}}_k\}} & \quad \sum_{l\in{\mc{N}}_k}a_{lk}^2\left[\mu_l^2
\sigma_{v,l}^2\Tr(R_{u,l})+\Tr(R_{v,lk}^{(\psi)})\right]\\
\st & \;\; \sum_{l\in{\mc{N}}_k}a_{lk}=1,\;\;a_{lk}\ge0,\;\;a_{lk}=0\;\mbox{if}\;l\notin{\mc{N}}_k  \\
\end{aligned}
}\end{equation}
for $k=1,\dots,N$. With each node $l\in\N_k$, we associate the following nonnegative \emph{variance product} measure:
\begin{align}
\label{eqn:gammalkdef}
\gamma_{lk}^2\defeq\begin{cases}
\mu_k^2\sigma_{v,k}^2\Tr(R_{u,k}), & l=k \\
\mu_l^2\sigma_{v,l}^2\Tr(R_{u,l})\!+\!{\mathrm{Tr}}(R_{v,lk}^{(\psi)}), & l\!\in\!{\mc{N}}_k\backslash\{k\} \\
\end{cases}
\end{align}
This measure incorporates information about the link noise covariances $\{R_{v,lk}^{(\psi)}\}$. The solution of \eqref{eqn:minsubproblem} is then given by
\begin{equation}
\label{eqn:relativevariance}
\boxed{
a_{lk}=\begin{cases}
\displaystyle\frac{\gamma_{lk}^{-2}}{\sum_{m\in{\mc{N}}_k}\gamma_{mk}^{-2}}, & {\textrm{if $l\in{\mc{N}}_k$}} \\
0, & {\textrm{otherwise}} \\
\end{cases}
} \quad \left(\begin{minipage}{1.3in}
relative variance rule
\end{minipage}\right)
\end{equation}
We refer to this combination rule as the relative variance combination rule; it is an extension of the rule devised in \cite{Tu11CAMSAP} to the case of noisy information exchanges. In particular, the definition of the scalars $\{\gamma_{lk}^2\}$ in \eqref{eqn:gammalkdef} is different and now depends on both subscripts $l$ and $k$.

Minimizing the EMSE expression \eqref{eqn:noisyatcEMSE} for the ATC algorithm over left-stochastic matrices $A$ can be pursued in a similar manner by noting that
\begin{align}
\label{eqn:emseatcup}
{\mathrm{Tr}}\left[{\mc{B}}_{\textrm{atc}}^j({\mc{A}}^{\T}{\mc{M}}{\mc{S}}{\mc{M}}
{\mc{A}}\!+\!{\mc{R}}_v^{(\psi)}){\mc{B}}_{\textrm{atc}}^{*j}{\mc{R}}_u\right]\le c^2[\rho(I_{NM}\!-\!{\mc{M}}{\mc{R}}_u)]^{2j}\,
{\mathrm{Tr}}({\mc{A}}^{\T}{\mc{M}}{\mc{S}}{\mc{M}}{\mc{A}}\!+\!{\mc{R}}_v^{(\psi)})\,{\mathrm{Tr}}({\mc{R}}_u)
\end{align}
Thus, minimizing the upper bound of the network EMSE leads to the same solution \eqref{eqn:relativevariance}. Using the same argument, we can also show that the same result minimizes the upper bound of the network MSD or EMSE for the CTA algorithm.

\subsection{Adaptive Combination Rule}
To apply the relative variance combination rule \eqref{eqn:relativevariance}, each node $k$ needs to know the variance products, $\{\gamma_{lk}^2\}$, of their neighbors, which in general are not available since they require knowledge of the quantities $\{\sigma_{v,l}^2,{\mathrm{Tr}}(R_{u,l}),{\mathrm{Tr}}(R_{v,lk}^{(\psi)})\}$. Therefore, we now propose an adaptive combination rule by using data that are available to the individual nodes. For the ATC algorithm, we first note from \eqref{eqn:noisypsilkidef} and \eqref{eqn:noisydiffusionincrementalold} that
\begin{align}
{\mathbb{E}}\|{\bs{\psi}}_{lk,i}-{\bs{w}}_{l,i-1}\|^2&\approx\mu_l^2\sigma_{v,l}^2\Tr(R_{u,l})+\Tr(R_{v,lk}^{(\psi)})\nonumber\\
{}&=\gamma_{lk}^2
\end{align}
for $l\in\N_k\backslash\{k\}$. Since the algorithm converges in the mean and mean-square senses under Assumption \ref{asm:smallstepsize}, all the estimates $\{{\bm{w}}_{k,i}\}$ tend close to $w^o$ as $i\rightarrow\infty$. This allows us to estimate $\gamma_{lk}^2$ for node $k$ by using instantaneous realizations of $\|{\bs{\psi}}_{lk,i}-{\bs{w}}_{k,i-1}\|^2$, where we replace ${\bs{w}}_{l,i-1}$ by ${\bs{w}}_{k,i-1}$. Similarly, for node $k$ itself, we can use realizations of $\|{\bs{\psi}}_{k,i}-{\bs{w}}_{k,i-1}\|^2$ to estimate $\gamma_{kk}^2$. To unify the notation, we define ${\bs{\psi}}_{kk,i}\defeq{\bs{\psi}}_{k,i}$. Let ${\widehat{\bm{\gamma}}}_{lk}^2(i)$ denote an estimator for $\gamma_{lk}^2$ that is computed by node $k$ at time $i$. Then, one way to evaluate ${\widehat{\bm{\gamma}}}_{lk}^2(i)$ is through the recursion:
\begin{equation}
\boxed{
\label{eqn:adaptivecombination}
{\widehat{\bm{\gamma}}}_{lk}^2(i)=(1-\nu_k){\widehat{\bm{\gamma}}}_{lk}^2(i-1)
+\nu_k\|{\bs{\psi}}_{lk,i}-{\bs{w}}_{k,i-1}\|^2
}\end{equation}
for $l\in{\mc{N}}_k$, where $\nu_k\in(0,1)$ is a forgetting factor that is usually close to one. In this way, we arrive at the adaptive combination rule:
\begin{equation}
\label{eqn:relativevarianceadapt}
\boxed{
\bm{a}_{lk}(i)=\begin{cases}
\displaystyle\frac{[{\widehat{\bm{\gamma}}}_{lk}^{2}(i)]^{-1}}{\sum_{m\in{\mc{N}}_k}[{\widehat{\bm{\gamma}}}_{mk}^{2}(i)]^{-1}}, & {\textrm{if $l\in{\mc{N}}_k$}} \\
0, & {\textrm{otherwise}} \\
\end{cases}
}\end{equation}

\section{Mean-Square Tracking Behavior}
The diffusion strategy \eqref{eqn:idealdiffusionpriordiff}--\eqref{eqn:idealdiffusionpostdiff} is adaptive in nature. One of the main benefits of adaptation (by using constant step-sizes) is that it endows networks with tracking abilities when the underlying weight vector $w^o$ varies with time. In this section we analyze how well an adaptive network is able to track variations in $w^o$. To do so, we adopt a random-walk model for $w^o$ that is commonly used in the literature to describe the non-stationarity of the weight vector \cite{Sayed08}.
\vspace{0.5\baselineskip}
\begin{assumption}[Random-walk model]
\label{asm:randomwalk} The weight vector $w^o$ changes according to the model:
\begin{align}
\label{eqn:randomwalk}
{\bs{w}}_{i}^o={\bs{w}}_{i-1}^o+{\bs{\eta}}_{i}
\end{align}
where $\{{\bs{w}}_{i}^o\}$ has a constant mean $w^o$ for all $i$, $\{{\bs{\eta}}_{i}\}$ is an i.i.d. random sequence with zero mean and covariance matrix $R_{\eta}$; the sequence $\{{\bs{\eta}}_{i}\}$ is independent of the initial conditions $\{{\bs{w}}_{-1}^o,{\bs{w}}_{k,-1}\}$ and of all regression data and noise signals across the network for all time instants. \hfill \IEEEQED
\end{assumption}
\vspace{0.5\baselineskip}
We now define the error vector at node $k$ as
\begin{align}
{\widetilde{\bs{w}}}_{k,i}\defeq{\bs{w}}_{i}^o-{\bs{w}}_{k,i}
\end{align}
so that the global error recursion (\ref{eqn:noisyerrorrecursion1}) for the network is replaced by
\begin{align}
\label{eqn:trackerrorrecursion1}
{\wt{\bs{w}}}_i&={\mc{A}}_2^{\T}(I_{NM}\!-\!{\mc{M}}{\bs{\mc{R}}}_i'){\mc{A}}_1^{\T}{\wt{\bs{w}}}_{i-1}
\!+\!{\mc{A}}_2^{\T}(I_{NM}\!-\!{\mc{M}}{\bs{\mc{R}}}_i'){\mc{A}}_1^{\T}{\bs{\zeta}}_i\nonumber\\
{}&\qquad-\!{\mc{A}}_2^{\T}(I_{NM}\!-\!{\mc{M}}{\bs{\mc{R}}}_i'){\bs{v}}_{i-1}^{(w)}
\!-\!{\mc{A}}_2^{\T}{\mc{M}}{\bs{z}}_i\!-\!{\bs{v}}_i^{(\psi)}
\end{align}
where the $NM\times1$ vector ${\bs{\zeta}}_i$ is defined as
\begin{align}
\label{eqn:zetadef}
{\bs{\zeta}}_i\defeq{\mathrm{col}}\left\{{\bs{\eta}}_{i},\dots,{\bs{\eta}}_{i}\right\}={\mathds{1}}_{N}\otimes{\bs{\eta}}_{i}
\end{align}

\subsection{Convergence Conditions}
By Assumptions \ref{asm:all} and \ref{asm:randomwalk}, it can be verified that the condition for mean convergence continues to be $\rho\left({\mc{B}}\right)<1$, where ${\mc{B}}$ is defined in \eqref{eqn:bigBdef}. In addition, it can also be verified that the error recursion \eqref{eqn:trackerrorrecursion1} converges in the mean sense to the same non-zero bias vector $g$ as in \eqref{eqn:gdef}. From \eqref{eqn:trackerrorrecursion1} and under Assumption \ref{asm:smallstepsize}, we can derive the weighted variance relation:
\begin{align}
\label{eqn:weightedvariancetrack}
{\mathbb{E}}\|{\wt{\bs{w}}}_i\|_{\sigma}^2&\approx{\mathbb{E}}\|{\wt{\bs{w}}}_{i-1}\|_{{\mc{F}}\sigma}^2
+{\mathbb{E}}\|{\mc{A}}_2^{\T}(I_{NM}\!-\!{\mc{M}}{\bs{\mc{R}}}_i'){\mc{A}}_1^{\T}{\bs{\zeta}}_i\|_{\sigma}^2\nonumber\\
{}&\qquad-2\,\Re\{\E[{\bs{z}}_i^*{\mc{M}}{\mc{A}}_2\Sigma{\mc{A}}_2^{\T}(I_{NM}\!-\!{\mc{M}}{\bs{\mc{R}}}_i'){\mc{A}}_1^{\T}{\wt{\bs{w}}}_{i-1}]\}\nonumber\\
{}&\qquad+{\mathbb{E}}\|{\mc{A}}_2^{\T}(I_{NM}-{\mc{M}}{\bs{\mc{R}}}_i'){\bs{v}}_{i-1}^{(w)}\|_{\sigma}^2\nonumber\\
{}&\qquad+{\mathbb{E}}\|{\mc{A}}_2^{\T}{\mc{M}}{\bs{z}}_i\|_{\sigma}^2+{\mathbb{E}}\|{\bs{v}}_i^{(\psi)}\|_{\sigma}^2
\end{align}
where ${\mc{F}}$ is given in \eqref{eqn:Fdefsmall}. If the step-sizes are sufficiently small, then we can assume that the network continues to be mean-square stable.

\subsection{Steady-State Performance}
The steady-state performance is affected by the non-stationarity of $w^o$. From Assumption \ref{asm:smallstepsize}, at steady-state, expression \eqref{eqn:weightedvariancetrack} becomes
\begin{align}
\lim_{i\rightarrow\infty}{\mathbb{E}}\|{\wt{\bs{w}}}_{i}\|_{\Omega}^2\!\approx\!
[\vec({\mc{A}}_2^{\T}{\mc{M}}{\mc{C}}^{\T}{\mc{S}}{\mc{C}}{\mc{M}}{\mc{A}}_2
\!+\!{\mc{A}}_2^{\T}{\mc{A}}_1^{\T}{\mc{R}}_\zeta{\mc{A}}_1{\mc{A}}_2
\!+\!{\mc{R}}_v\!+\!{\mc{Y}}\!+\!{\mc{Y}}^*)]^*(I_{N^2M^2}\!-\!{\mc{F}})^{-1}{\mathrm{vec}}(\Omega)
\end{align}
where ${\mc{S}}$ is given in \eqref{eqn:bigSdef}, $R_v$ in \eqref{eqn:bigRvdef}, ${\mc{Y}}$ in \eqref{eqn:bigYdef}, ${\mc{F}}$ in \eqref{eqn:Fdefsmall}, and ${\mc{R}}_{\zeta}$ is the covariance matrix of ${\bs{\zeta}}_i$:
\begin{align}
\label{eqn:Rzetadef}
{\mc{R}}_{\zeta}\defeq{\mathbb{E}}{\bs{\zeta}}_i{\bs{\zeta}}_i^*=({\mathds{1}}_{N}{\mathds{1}}_{N}^{\T})\otimes R_{\eta}
\end{align}
By \eqref{eqn:A1A2def}, \eqref{eqn:bigA1CA2def}, and \eqref{eqn:Rzetadef}, we get
\begin{align}
\label{eqn:Retaiden}
{\mc{A}}_2^{\T}{\mc{A}}_1^{\T}{\mc{R}}_\zeta{\mc{A}}_1{\mc{A}}_2
&=(A_2^{\T}A_1^{\T}{\mathds{1}}_{N}{\mathds{1}}_{N}^{\T}A_1A_2)\otimes R_{\eta}\nonumber\\
{}&=({\mathds{1}}_{N}{\mathds{1}}_{N}^{\T})\otimes R_{\eta}\nonumber\\
{}&={\mc{R}}_{\zeta}
\end{align}
Then, following the same argument that led to \eqref{eqn:noisyMSD}, we find that the network MSD is now given by:
\begin{equation}
\boxed{
\label{eqn:trackMSD}
\begin{aligned}
{\textrm{MSD}}_{\textrm{trk}}&\!\approx\!\frac{1}{N}\![\vec({\mc{A}}_2^{\T}{\mc{M}}{\mc{C}}^{\T}{\mc{S}}{\mc{C}}{\mc{M}}{\mc{A}}_2
\!+\!{\mc{R}}_\zeta\!+\!{\mc{R}}_v\!+\!{\mc{Y}}\!+\!{\mc{Y}}^*)]^*(I_{N^2M^2}-{\mc{F}})^{-1}\vec(I_{NM})
\end{aligned}
}\end{equation}
Similarly, the network EMSE is given by:
\begin{equation}
\boxed{
\label{eqn:trackEMSE}
\begin{aligned}
{\textrm{EMSE}}_{\textrm{trk}}&\!\approx\!\frac{1}{N}\![\vec({\mc{A}}_2^{\T}{\mc{M}}{\mc{C}}^{\T}{\mc{S}}{\mc{C}}{\mc{M}}{\mc{A}}_2
\!+\!{\mc{R}}_\zeta\!+\!{\mc{R}}_v\!+\!{\mc{Y}}\!+\!{\mc{Y}}^*)]^*(I_{N^2M^2}-{\mc{F}})^{-1}\vec({\mc{R}}_u)
\end{aligned}
}\end{equation}
where ${\mc{R}}_u$ is defined in \eqref{eqn:bigRudef}. Observe that the main difference relative to \eqref{eqn:noisyMSD} and \eqref{eqn:noisyEMSE} is the addition of the term ${\mc{R}}_{\zeta}$. Therefore, all the results that were derived in the earlier section, such as \eqref{eqn:noisyMSD1} and \eqref{eqn:noisyEMSE1}, continue to hold by adding ${\mc{R}}_{\zeta}$. In particular, if Assumptions \ref{asm:smallstepsize} and \ref{asm:nodatasharing} are adopted, expressions \eqref{eqn:trackMSD} and \eqref{eqn:trackEMSE} can be approximated as
\begin{equation}
\boxed{
\label{eqn:trackMSD2}
\begin{aligned}
{\textrm{MSD}}_{\textrm{trk}}&\approx\frac{1}{N}[\vec({\mc{A}}_2^{\T}{\mc{M}}{\mc{S}}{\mc{M}}{\mc{A}}_2+{\mc{R}}_\zeta
+{\mc{R}}_v)]^*(I_{N^2M^2}-{\mc{F}})^{-1}\vec(I_{NM})
\end{aligned}
}\end{equation}
and
\begin{equation}
\boxed{
\label{eqn:trackEMSE2}
\begin{aligned}
{\textrm{EMSE}}_{\textrm{trk}}&\approx\frac{1}{N}[\vec({\mc{A}}_2^{\T}{\mc{M}}{\mc{S}}{\mc{M}}{\mc{A}}_2+{\mc{R}}_\zeta
+{\mc{R}}_v)]^*(I_{N^2M^2}-{\mc{F}})^{-1}\vec({\mc{R}}_u)
\end{aligned}
}\end{equation}
where ${\mc{R}}_v$ is now given in \eqref{eqn:bigRvsimple}.

\section{Simulation Results}
We simulate two scenarios: noisy information exchanges and non-stationary environments. We consider a connected network with $N=20$ nodes. The network topology is shown in Fig. \ref{fig:topology}.

\subsection{Imperfect Information Exchange}
\begin{figure}[t]
\centering
\includegraphics[height=1.8in]{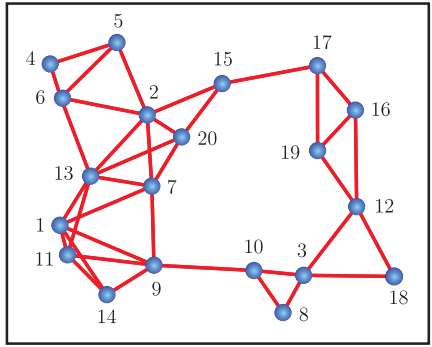}
\caption{A network topology with $N=20$ nodes.}
\label{fig:topology}
\end{figure}

\begin{figure}[t]
\centerline{
\subfloat[The variance profile of regression data.]
{\includegraphics[height=1.2in]{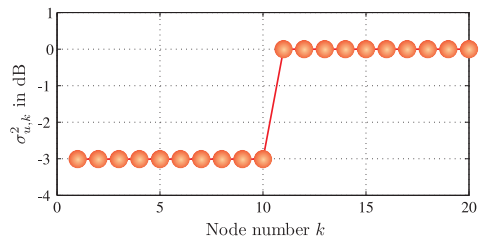}
\label{fig:regressor}}}
\centerline{
\subfloat[The variance profile of measurement noises.]
{\includegraphics[height=1.2in]{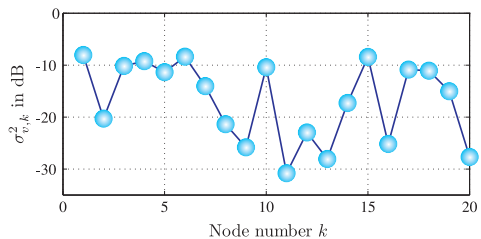}
\label{fig:modelnoise}}}
\caption{The variance profiles for regression data and measurement noises.}
\label{fig:sim1}
\end{figure}

\begin{figure*}[t!]
\centering
\includegraphics[height=3.5in]{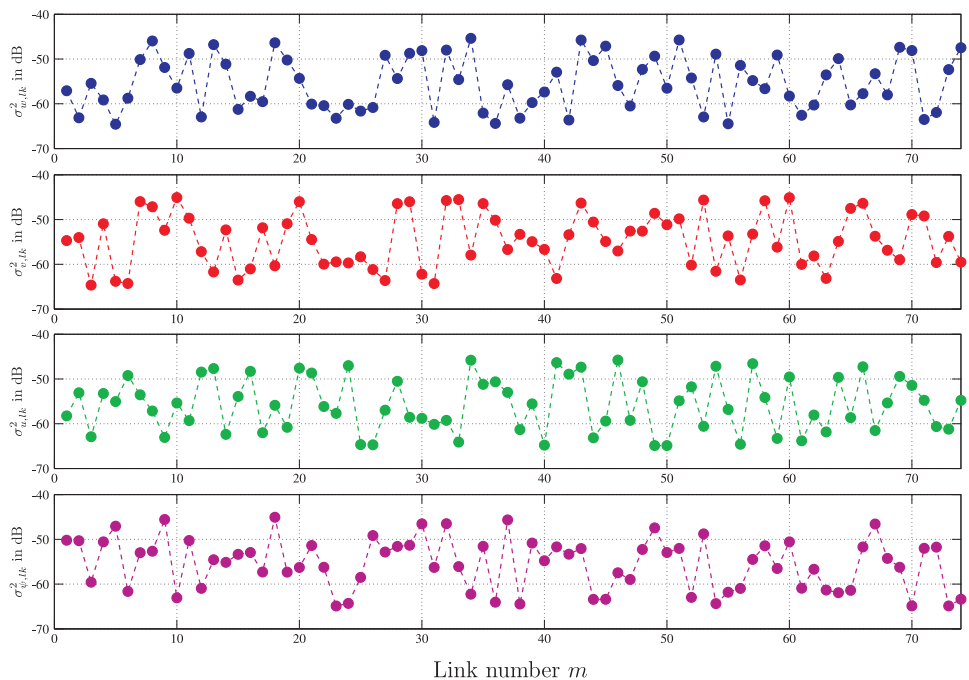}
\caption{The variance profiles for various sources of link noises, including $\{\sigma_{w,lk}^2,\sigma_{v,lk}^{2},\sigma_{u,lk}^2,\sigma_{\psi,lk}^2\}$.}
\label{fig:linknoise}
\vspace{0\baselineskip}
\end{figure*}

The unknown complex parameter $w^o$ of length $M=2$ is randomly generated; its value is $[0.3750+j2.0834,0.7174+j1.4123]$. We adopt uniform step-sizes, $\{\mu_k=0.01\}$, and uniformly white Gaussian regression data with covariance matrices $\{R_{u,k}=\sigma_{u,k}^2I_M\}$, where $\{\sigma_{u,k}^2\}$ are shown in Fig. \ref{fig:regressor}. The variances of the model noises, $\{\sigma_{v,k}^{2}\}$, are randomly generated and shown in Fig. \ref{fig:modelnoise}. We also use white Gaussian link noise signals such that $R_{v,lk}^{(w)}=\sigma_{w,lk}^2I_M$, $R_{v,lk}^{(u)}=\sigma_{u,lk}^2I_M$, and $R_{v,lk}^{(\psi)}=\sigma_{\psi,lk}^2I_M$. All link noise variances, $\{\sigma_{w,lk}^2,\sigma_{v,lk}^{2},\sigma_{u,lk}^2,\sigma_{\psi,lk}^2\}$, are randomly generated and illustrated in Fig. \ref{fig:linknoise} from top to bottom. We assign the link number by the following procedure. We denote the link from node $l$ to node $k$ as $\ell_{l,k}$, where $l\ne k$. Then, we collect the links $\{\ell_{l,k},l\in\N_k\backslash\{k\}\}$ in an ascending order of $l$ in the list ${\mc{L}}_k$ (which is a set with \emph{ordered} elements) for each node $k$. For example, for node $k=2$ in Fig. \ref{fig:topology}, it has $6$ links; the ordered links are then collected in ${\mc{L}}_2\defeq\{\ell_{5,2},\ell_{6,2},\ell_{7,2},\ell_{13,2},\ell_{15,2},\ell_{20,2}\}$. We concatenate $\{{\mc{L}}_k\}$ in an ascending order of $k$ to get the overall list ${\mc{L}}\defeq\{{\mc{L}}_1,{\mc{L}}_2,\dots,{\mc{L}}_N\}$. Eventually, the $m$th link in the network is given by the $m$th element in the list ${\mc{L}}$.

\begin{figure*}[t!]
\centerline{
\subfloat[Network MSD curves for ATC algorithms]
{\includegraphics[width=3in]{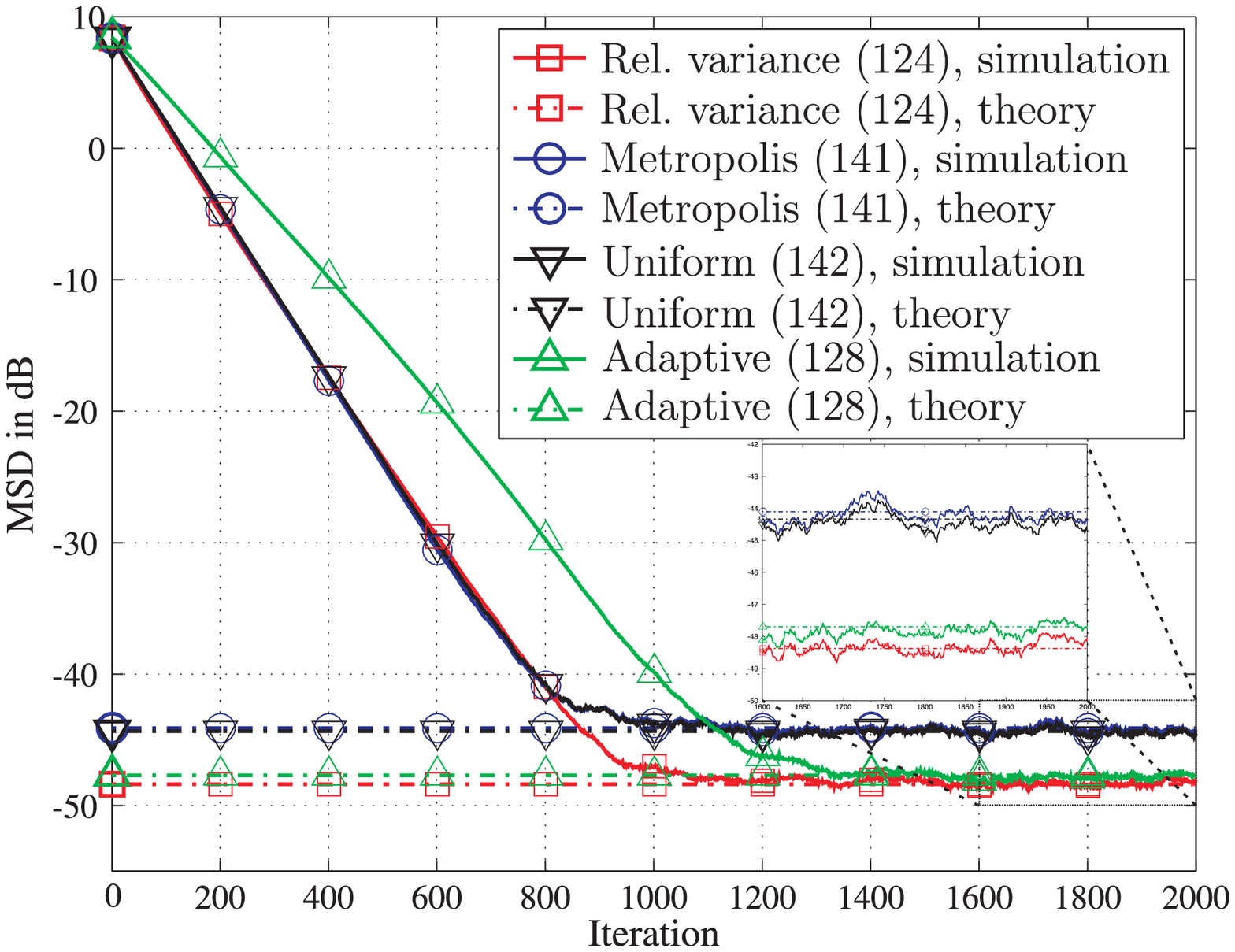}
\label{fig:msd_atc}}
\hfil
\subfloat[Network MSD curves for CTA algorithms]
{\includegraphics[width=3in]{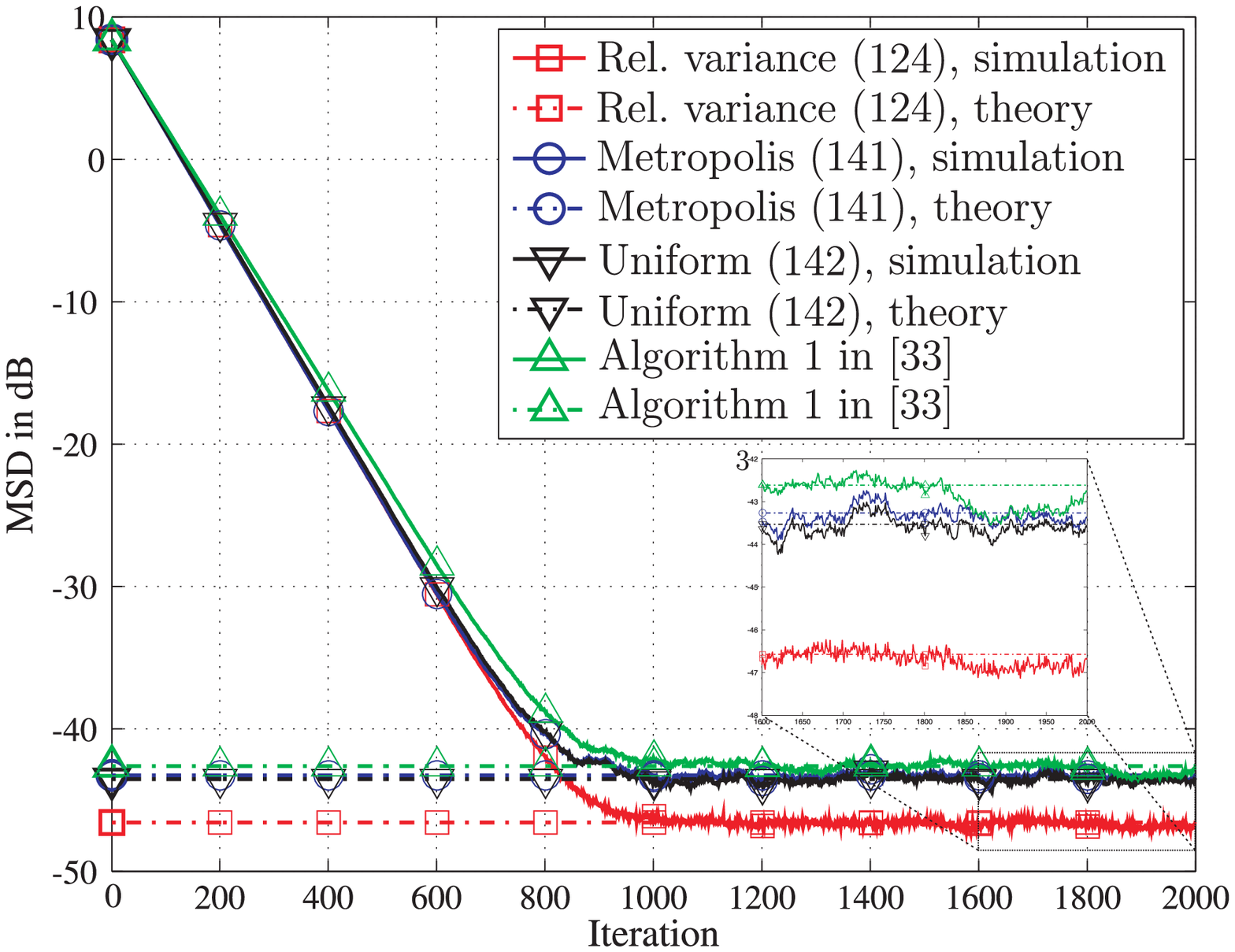}
\label{fig:msd_cta}}}
\centerline{
\subfloat[Network EMSE curves for ATC algorithms]
{\includegraphics[width=3in]{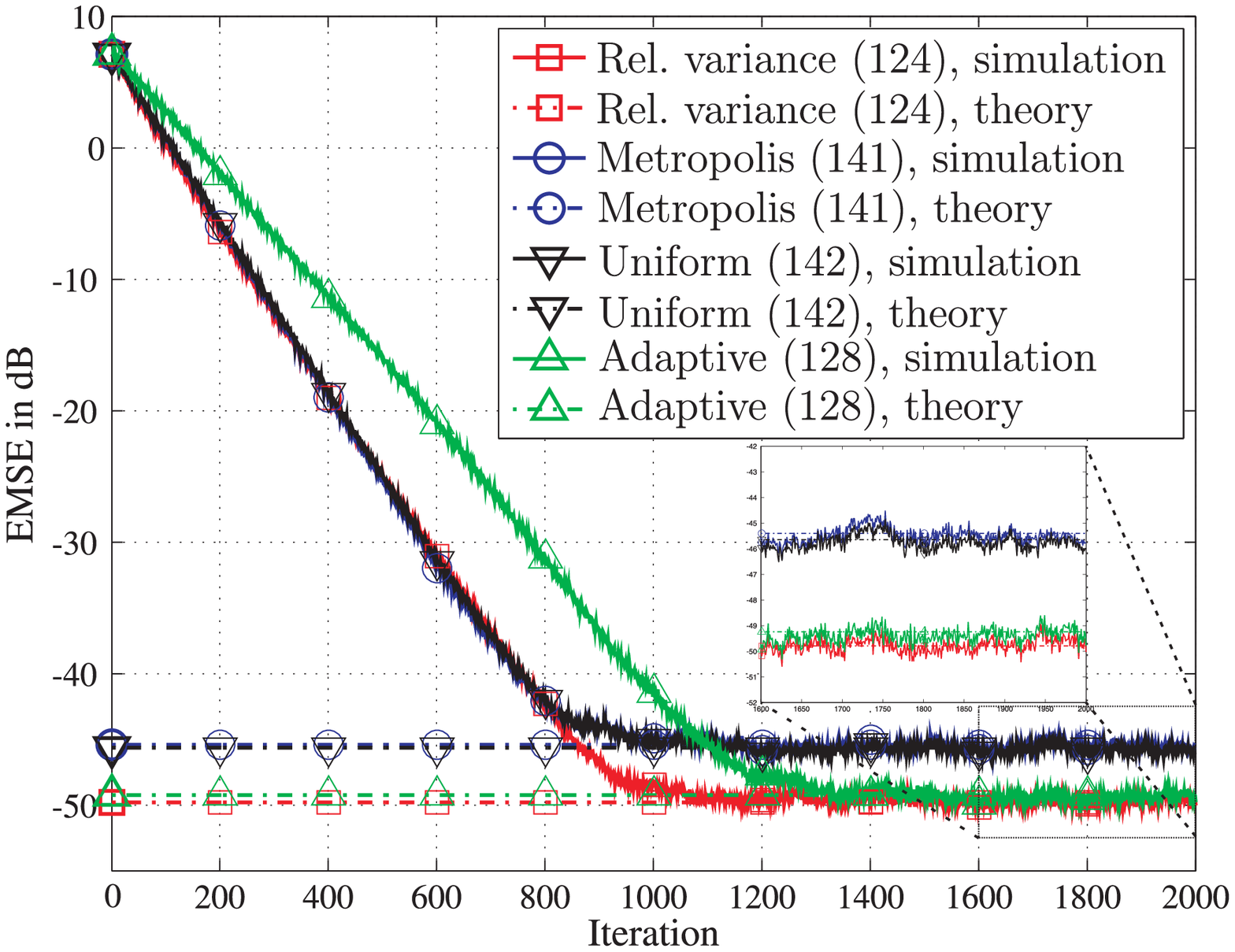}
\label{fig:emse_atc}}
\hfil
\subfloat[Network EMSE curves for CTA algorithms]
{\includegraphics[width=3in]{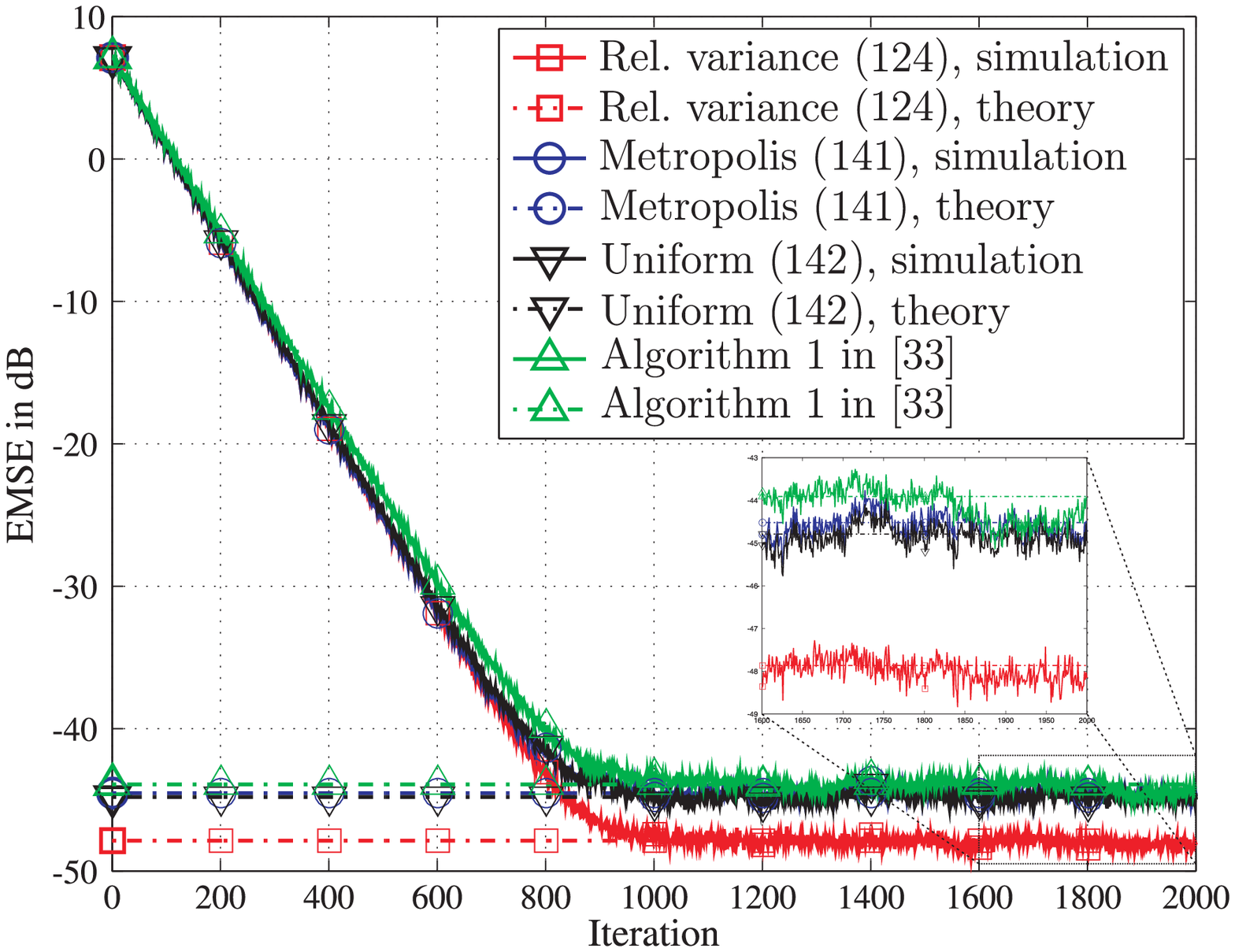}
\label{fig:emse_cta}}}
\caption{Simulated network MSD and EMSE curves and theoretical results (\ref{eqn:noisyMSD1}) and (\ref{eqn:noisyEMSE1}) for diffusion algorithms with various combination rules under noisy information exchange.}
\label{fig:sim2}
\vspace{-1\baselineskip}
\end{figure*}

We examine the simplified CTA and ATC algorithms in \eqref{eqn:ctaalgorithm} and \eqref{eqn:atcalgorithm}, namely, no sharing of data among nodes (i.e., $C=I_N$), under various combination rules: (i) the relative variance rule in \eqref{eqn:relativevariance}, (ii) the Metropolis rule in \cite{Metropolis53JCP}:
\begin{align}
\label{eqn:reldegvar}
\begin{cases}
\displaystyle
    a_{lk}=\frac{1}{\max\{|{\mc{N}}_{k}|,|{\mc{N}}_{l}|\}},&l\in{\mc{N}}_k\backslash\{k\}\\
\displaystyle    a_{lk}=1-\sum_{l\in{\mc{N}}_{k}\backslash\{k\}}a_{lk}, & l=k \\
    a_{lk}=0,&l\notin{\mc{N}}_k
  \end{cases}
\end{align}
where $|{\mc{N}}_{k}|$ denotes the degree of node $k$ (including the node itself), (iii) the uniform weighting rule:
\begin{align}
\label{eqn:uniformweights}
\begin{cases}
\displaystyle    a_{lk}=\frac{1}{|{\mc{N}}_k|},&l\in{\mc{N}}_k\\
    a_{lk}=0,&l\notin{\mc{N}}_k
  \end{cases}
\end{align}
and (iv) the adaptive rule in \eqref{eqn:relativevarianceadapt} with $\{\nu_k=0.05\}$. We plot the network MSD and EMSE learning curves for ATC algorithms in Figs. \ref{fig:msd_atc} and \ref{fig:emse_atc} by averaging over 50 experiments. For CTA algorithms, we plot their network MSD and EMSE learning curves in Figs. \ref{fig:msd_cta} and \ref{fig:emse_cta} also by averaging over 50 experiments. Moreover, we also plot their theoretical results \eqref{eqn:noisyMSD1} and \eqref{eqn:noisyEMSE1} in the same figures. From Fig. \ref{fig:sim2} we see that the relative variance rule makes diffusion algorithms achieve the lowest MSD and EMSE levels at steady-state, compared to the metropolis and uniform rules as well as the algorithm from \cite{Abdolee11DCOSS} (which also requires knowledge of the noise variances). In addition, the adaptive rule attains MSD and EMSE levels that are only slightly larger than those of the relative variance rule, although, as expected, it converges slower due to the additional learning step \eqref{eqn:adaptivecombination}.

\subsection{Non-stationary Scenario}
\begin{figure}[t]
\centerline{
\subfloat[The low noise-level case.]
{\includegraphics[height=2.5in]{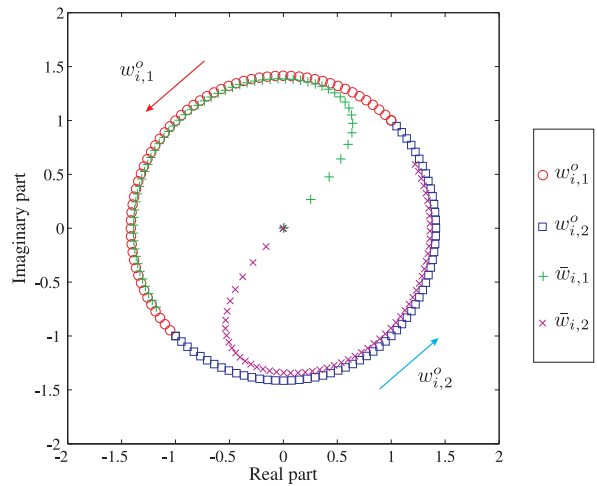}
\label{fig:track_hsnr}}}
\centerline{
\subfloat[The high noise-level case.]
{\includegraphics[height=2.5in]{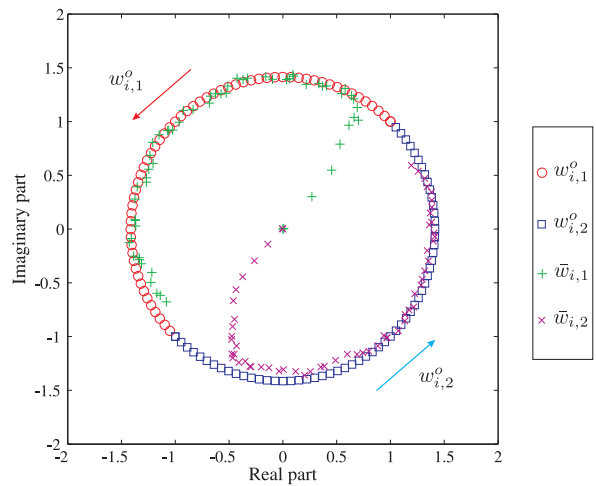}
\label{fig:track_lsnr}}}
\caption{An adaptive network tracking a parameter vector $w^o\in{\mathbb{C}}^{2}$.}
\label{fig:track}
\end{figure}

\begin{figure}[t]
\centerline{
\subfloat[The variance profile for low noise-level.]
{\includegraphics[height=1.2in]{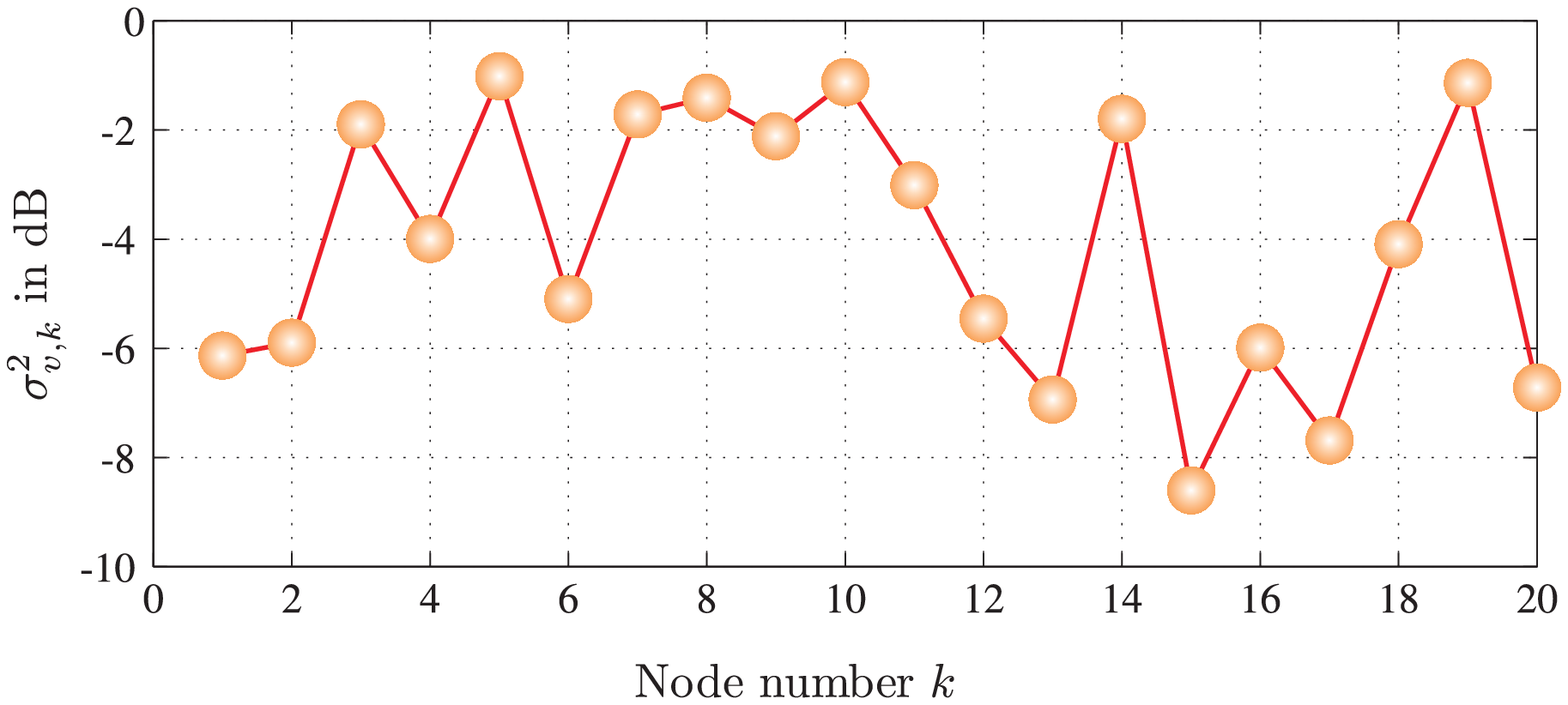}
\label{fig:hsnr}}}
\centerline{
\subfloat[The variance profile for high noise-level.]
{\includegraphics[height=1.2in]{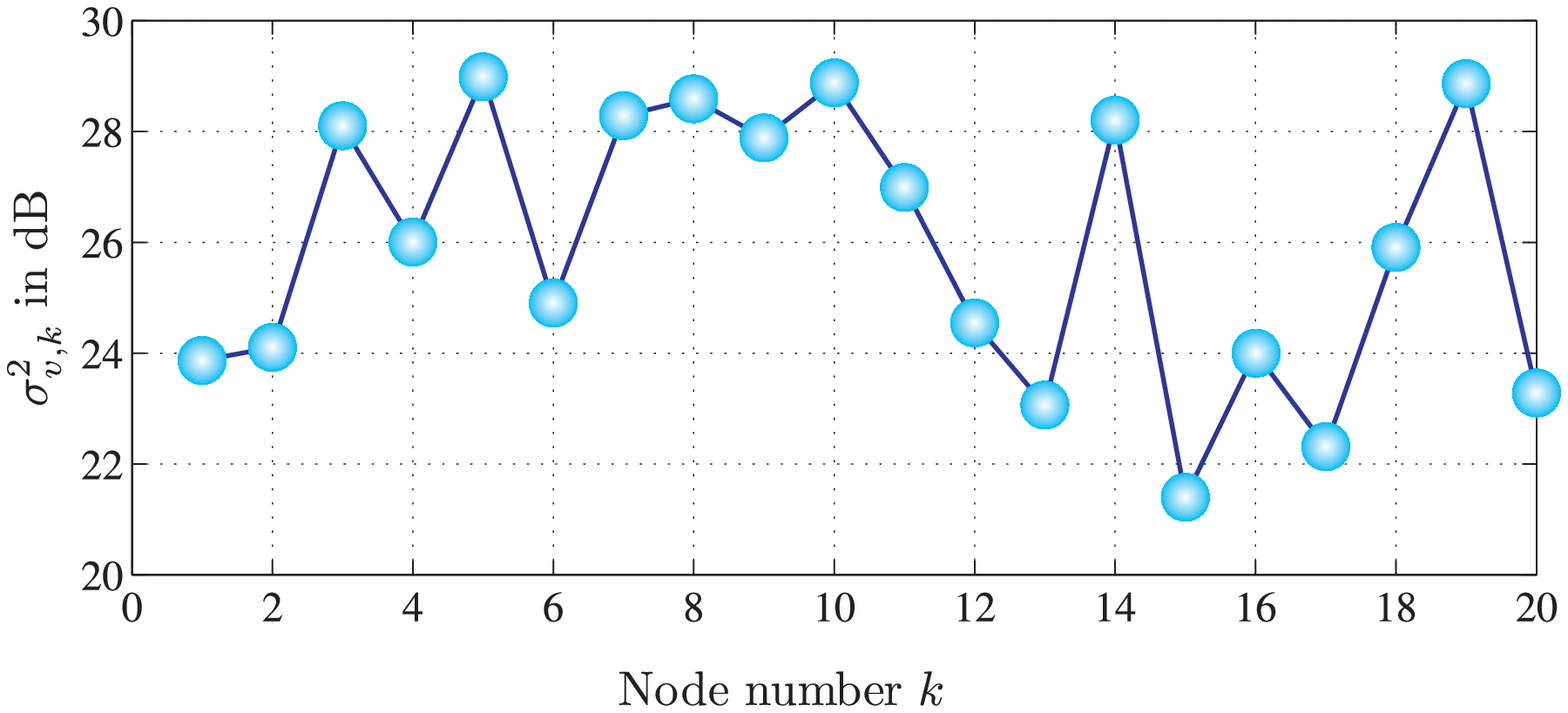}
\label{fig:lsnr}}}
\caption{The noise variance profiles for two cases.}
\label{fig:sim3}
\end{figure}

The value for each entry of the complex parameter $w_i^o={\mathrm{col}}\{w_{i,1}^o,w_{i,2}^o\}$ is assumed to be changing over time along a circular trajectory in the complex plane, as shown in Fig. \ref{fig:track}. The dynamic model for $w_i^o$ is expressed as $w_{i,m}^o=e^{j\omega}w_{i-1,m}^o$, where $m=1,2$, $\omega=2\pi/6000$, and $w_{-1}^o=\col\{1+j,-1-j\}$. The covariance matrices $\{R_{u,k}\}$ are randomly generated such that $R_{u,k}\ne R_{u,l}$ when $k\ne l$, but their traces are normalized to be one, i.e., $\Tr(R_{u,k})=1$, for all nodes. The variances for the model noises, $\{\sigma_{v,k}^2\}$, are also randomly generated. We examine two different scenarios: the low noise-level case where the average noise variance across the network is $-5$ dB and the noise variances are shown in Fig. \ref{fig:hsnr}; and the high noise-level case where the average variance is $25$ dB and the variances are shown in Fig. \ref{fig:lsnr}. We simulate 3000 iterations and average over 20 experiments in Figs. \ref{fig:track_hsnr} and \ref{fig:track_lsnr} for each case. The step-size is 0.01 and uniform across the network. For simplicity, we adopt the simplified ATC algorithm where $C=I_N$, and only use the uniform weighting rule \eqref{eqn:uniformweights}. The tracking behavior of the network, denoted as ${\bar{w}}_i={\mathrm{col}}\{{\bar{w}}_{i,1},{\bar{w}}_{i,2}\}$, is obtained by averaging over all the estimates, $\{w_{k,i}\}$, across the network. Figs. \ref{fig:track_hsnr} and \ref{fig:track_lsnr} depict the complex plane; the horizontal axis is the real axis and the vertical axis is the imaginary axis. Therefore, for every time $i$, each entry of $w_i^o$ or ${\bar{w}}_i$ represents a point in the plane. When $i$ is increasing, $w_{i,1}^o$ moves along the red trajectory (in $\circ$), $w_{i,2}^o$ along the blue trajectory (in $\square$), ${\bar{w}}_{i,1}$ along the green trajectory (in $+$), and ${\bar{w}}_{i,2}$ along the magenta trajectory (in $\times$). From Fig. \ref{fig:track}, it can be seen that diffusion algorithms exhibit the tracking ability in both high and low noise-level environments.

\section{Conclusions}
In this work we investigated the performance of diffusion algorithms under several sources of noise during information exchange and under non-stationary environments. We first showed that, on one hand, the link noise over the regression data biases the estimators and deteriorates the conditions for mean and mean-square convergence. On the other hand, diffusion strategies can still stabilize the mean and mean-square convergence of the network with noisy information exchange. We derived analytical expressions for the network MSD and EMSE and used these expressions to motivate the choice of combination weights that help ameliorate the effect of information-exchange noise and improve network performance. We also extended the results to the non-stationary scenario where the unknown parameter $w^o$ is changing over time. Simulation results illustrate the theoretical findings and how well they match with theory.

\appendices

\section{Stability of ${\mc{A}}_2^{\T}\left(I_{NM}-{\mc{M}}{\mc{R}}'\right){\mc{A}}_1^{\T}$}
\label{app:meanconvergence}
Following \cite{Takahashi10TSP}, we first define the block maximum norm of a vector.
\vspace{0.5\baselineskip}
\begin{definition}[Block Maximum Norm]
\label{def:vecblkmaxnorm}
Given a vector $x={\mathrm{col}}\{x_1,\dots,x_N\}\in{\mathbb{C}}^{MN}$ consisting of $N$ blocks $\{x_k\in{\mathbb{C}}^{M},k=1,\dots,N\}$, the block maximum norm is the real function $\|\cdot\|_{b,\infty}:\;{\mathbb{C}}^{MN}\rightarrow{\mathbb{R}}$, defined as
\begin{align}
\label{eqn:blkmaxnormdef}
\|x\|_{b,\infty}\defeq\max_{1\le k\le N}\|x_k\|_2
\end{align}
where $\|\cdot\|_2$ denotes the standard $2$-norm on ${\mathbb{C}}^{M}$. \hfill \IEEEQED
\end{definition}
\vspace{0.5\baselineskip}
Similarly, we define the matrix norm that is induced by the block maximum norm as follows:
\vspace{0.5\baselineskip}
\begin{definition}[Block Maximum Matrix Norm]
\label{def:blkmaxnorm}
Given a block matrix ${\mc{A}}\in{\mathbb{C}}^{MN\times MN}$ with block size $M\times M$, then
\begin{align}
\label{eqn:matblkmaxnormdef}
\|{\mc{A}}\|_{b,\infty}\defeq\max_{x\in{\mathbb{C}}^{MN}\backslash\{0\}}\frac{\|{\mc{A}}x\|_{b,\infty}}{\|x\|_{b,\infty}}
\end{align}
denotes the induced block maximum (matrix) norm on ${\mathbb{C}}^{MN\times MN}$. \hfill \IEEEQED
\end{definition}
\vspace{0.5\baselineskip}
\begin{lemma}
\label{lemma:blkunitaryinvariant}
The block maximum matrix norm is block unitary invariant, i.e., given a block diagonal unitary matrix ${\mc{U}}\defeq{\mathrm{diag}}\{U_1,\dots,U_N\}\in{\mathbb{C}}^{MN\times MN}$ consisting of $N$ unitary blocks $\{U_k\in{\mathbb{C}}^{M\times M},k=1,\dots,N\}$, where $U_kU_k^*=U_k^*U_k=I_M$, for any matrix ${\mc{A}}\in{\mathbb{C}}^{MN\times MN}$, then
\begin{align}
\label{eqn:blkunitaryinv}
\|{\mc{A}}\|_{b,\infty}=\|\,{\mc{U}}{\mc{A}}\,{\mc{U}}^*\|_{b,\infty}
\end{align}
where $\|\cdot\|_{b,\infty}$ denotes the block maximum matrix norm on ${\mathbb{C}}^{MN\times MN}$ with block size $M\times M$.\hfill \IEEEQED
\end{lemma}
\vspace{0.5\baselineskip}
\begin{lemma}
\label{lemma:rghtstocmat}
Let $A\in{\mathbb{C}}^{N\times N}$ be a right-stochastic matrix. Then, for block size $M\times M$,
\begin{align}
\|A\otimes I_M\|_{b,\infty}=1
\end{align}
\end{lemma}
\begin{IEEEproof}
From Definition \ref{def:blkmaxnorm}, we get
\begin{align}
\label{eqn:rghtstocmatub}
\|A\otimes I_M\|_{b,\infty}&=\max_{x\in{\mathbb{C}}^{MN}\backslash\{0\}}
\frac{\max_{l}\|\sum_{k=1}^{N}[A]_{lk}x_k\|_2}{\max_{k}\|x_k\|_2}\nonumber\\
{}&\le\max_{x\in{\mathbb{C}}^{MN}\backslash\{0\}}
\frac{\max_{l}\sum_{k=1}^{N}[A]_{lk}\|x_k\|_2}{\max_{k}\|x_k\|_2}\nonumber\\
{}&\le\max_{x\in{\mathbb{C}}^{MN}\backslash\{0\}}
\frac{\max_{l}(\sum_{k=1}^{N}[A]_{lk})\cdot\max_k\|x_k\|_2}{\max_{k}\|x_k\|_2}\nonumber\\
{}&\le\max_{x\in{\mathbb{C}}^{MN}\backslash\{0\}}\frac{\max_{l}1\cdot\max_{k}\|x_k\|_2}{\max_{k}\|x_k\|_2}\nonumber\\
{}&=1
\end{align}
where $x\defeq{\mathrm{col}}\{x_1,\dots,x_N\}\in{\mathbb{C}}^{MN}$ consists of $N$ blocks $\{x_k\in{\mathbb{C}}^{M},k=1,\dots,N\}$, and $[A]_{lk}$ denotes the $(l,k)$th entry of $A$. On the other hand, for any induced matrix norm, say, the block maximum norm, it is always lower bounded by the spectral radius of the matrix \cite{Horn85}:
\begin{align}
\label{eqn:rghtstocmatlb}
\|A\otimes I_M\|_{b,\infty}\ge\rho(A\otimes I_M)=\rho(A)=1
\end{align}
Combining \eqref{eqn:rghtstocmatub} and \eqref{eqn:rghtstocmatlb} completes the proof.
\end{IEEEproof}
\vspace{0.5\baselineskip}
\begin{lemma}
\label{lemma:blockdiagonal}
Let ${\mc{A}}\in{\mathbb{C}}^{NM\times NM}$ be a block diagonal Hermitian matrix with block size $M\times M$. Then the block maximum norm of the matrix ${\mc{A}}$  is equal to its spectral radius, i.e.,
\begin{align}
\|{\mc{A}}\|_{b,\infty}=\rho({\mc{A}})
\end{align}
\end{lemma}
\begin{IEEEproof}
Denote the $k$th $M\times M$ submatrix on the diagonal of ${\mc{A}}$ by $A_k$. Let $A_k=U_k\Lambda_kU_k^*$ be the eigen-decomposition of $A_k$, where $U_k\in{\mathbb{C}}^{M\times M}$ is unitary and $\Lambda_k\in{\mathbb{R}}^{M\times M}$ is diagonal. Define the block unitary matrix ${\mc{U}}\defeq{\mathrm{diag}}\{U_1,\dots,U_N\}$ and the diagonal matrix $\Lambda\defeq{\mathrm{diag}}\{\Lambda_1,\dots,\Lambda_N\}$. Then, ${\mc{A}}={\mc{U}}\Lambda{\mc{U}}^*$. By Lemma \ref{lemma:blkunitaryinvariant}, the block maximum norm of ${\mc{A}}$ with block size $M\times M$ is
\begin{align}
\label{eqn:blockdiagonalA}
\|{\mc{A}}\|_{b,\infty}&=\|{\mc{U}}\Lambda{\mc{U}}^*\|_{b,\infty}\nonumber\\
{}&=\|\Lambda\|_{b,\infty}\nonumber\\
{}&=\max_{x\in{\mathbb{C}}^{MN}\backslash\{0\}}\frac{\max_{k}\|\Lambda_kx_k\|_2}{\max_{k}\|x_k\|_2}\nonumber\\
{}&\le\max_{x\in{\mathbb{C}}^{MN}\backslash\{0\}}\frac{\max_{k}\|\Lambda_k\|_2\cdot\|x_k\|_2}{\max_{k}\|x_k\|_2}\nonumber\\
{}&\le\max_{x\in{\mathbb{C}}^{MN}\backslash\{0\}}\frac{\max_{k}\|\Lambda_k\|_2\cdot\max_{k}\|x_k\|_2}{\max_{k}\|x_k\|_2}\nonumber\\
{}&=\max_{k}\|\Lambda_k\|_2\nonumber\\
{}&=\rho({\mc{A}})
\end{align}
where we used the fact that the induced $2$-norm is identical to the spectral radius for Hermitian matrices \cite{Horn85}. On the other hand, any matrix norm is lower bounded by the spectral radius \cite{Horn85}, i.e.,
\begin{align}
\label{eqn:blockdiagonalA1}
\rho({\mc{A}})\le\|{\mc{A}}\|_{b,\infty}
\end{align}
Combining \eqref{eqn:blockdiagonalA} and \eqref{eqn:blockdiagonalA1} completes the proof.
\end{IEEEproof}
\vspace{0.5\baselineskip}
Now we show that the matrix ${\mc{A}}_2^{\T}(I_{NM}-{\mc{M}}{\mc{R}}'){\mc{A}}_1^{\T}$ is stable if $I_{NM}-{\mc{M}}{\mc{R}}'$ is stable. For any induced matrix norm, say, the block maximum norm with block size $M\times M$, we have \cite{Horn85}
\begin{align}
\label{eqn:rhoAbounded}
\rho\left({\mc{A}}_2^{\T}(I_{NM}-{\mc{M}}{\mc{R}}'){\mc{A}}_1^{\T}\right)&
\le\|{\mc{A}}_2^{\T}(I_{NM}-{\mc{M}}{\mc{R}}'){\mc{A}}_1^{\T}\|_{b,\infty}\nonumber\\
{}&\le\|{\mc{A}}_2^{\T}\|_{b,\infty}\cdot\|I_{NM}-{\mc{M}}{\mc{R}}'\|_{b,\infty}\cdot\|{\mc{A}}_1^{\T}\|_{b,\infty}\nonumber\\
{}&=\|I_{NM}-{\mc{M}}{\mc{R}}'\|_{b,\infty}
\end{align}
where, from \eqref{eqn:A1A2def} and \eqref{eqn:bigA1CA2def}, ${\mc{A}}_1^{\T}$ and ${\mc{A}}_2^{\T}$ satisfy Lemma \ref{lemma:rghtstocmat}. By \eqref{eqn:bigMdef} and \eqref{eqn:bigRprimedef}, it is straightforward to see that $I_{NM}-{\mc{M}}{\mc{R}}'$ is block diagonal with block size $M\times M$. Then, by Lemma \ref{lemma:blockdiagonal}, expression \eqref{eqn:rhoAbounded} can be further expressed as
\begin{align}
\rho\left({\mc{A}}_2^{\T}(I_{NM}-{\mc{M}}{\mc{R}}'){\mc{A}}_1^{\T}\right)
\le\rho(I_{NM}-{\mc{M}}{\mc{R}}')
\end{align}
which completes the proof.\footnote{This statement fixes the argument that appeared in Appendix I of \cite{Cattivelli10TSP} and Lemma 2 of \cite{Cattivelli10TAC}. Since the matrix ${\mc{X}}$ in Appendix I of \cite{Cattivelli10TSP} and the matrix ${\mc{M}}$ in Lemma 2 of \cite{Cattivelli10TAC} are block diagonal, the $\|\cdot\|_{\rho}$ norm used in these references should simply be replaced by the $\|\cdot\|_{b,\infty}$ norm used here and as already done in \cite{Takahashi10TSP}.}

\section{Proof of expression \eqref{eqn:bigRzdef}}
\label{app:Rz}
Let us denote the $(l,k)$th submatrix of $R_{z}$ by $R_{z,lk}\in{\mathbb{C}}^{M\times M}$. By Assumptions \ref{asm:all} and expression \eqref{eqn:zkidef}, $R_{z,lk}$ can be evaluated as
\begin{align}
\label{eqn:appCovZterm1}
R_{z,lk}&=\E\,{\bs{z}}_{l,i}{\bs{z}}_{k,i}^*\nonumber\\
{}&=\sum_{m\in\N_l}\sum_{n\in\N_k}c_{ml}c_{nk}
\underbrace{\E\left({\bs{u}}_{ml,i}^*{\bs{v}}_{ml}(i){\bs{v}}_{nk}^*(i){\bs{u}}_{nk,i}\right)}_{\defeq\;R_{ml,nk}}
\end{align}
where, by expressions \eqref{eqn:noisyulkidef} and \eqref{eqn:vlkidef},
\begin{align}
\label{eqn:appCovZterm2}
R_{ml,nk}&={\mathbb{E}}\!\left({\bs{u}}_{m,i}\!+\!{\bs{v}}_{ml,i}^{(u)}\right)^*
\left({\bs{v}}_{m}(i)\!+\!{\bs{v}}_{ml}^{(d)}(i)\!-\!{\bs{v}}_{ml,i}^{(u)}w^o\right)
\left({\bs{v}}_{n}(i)\!+\!{\bs{v}}_{nk}^{(d)}(i)\!-\!{\bs{v}}_{nk,i}^{(u)}w^o\right)^*
\left({\bs{u}}_{n,i}\!+\!{\bs{v}}_{nk,i}^{(u)}\right)
\end{align}
When $m\ne n$, expression \eqref{eqn:appCovZterm2} reduces to
\begin{align}
\label{eqn:appCovZterm3}
R_{ml,nk}=R_{v,ml}^{(u)}w^ow^{o*}R_{v,nk}^{(u)}
\end{align}
When $m=n$, expression \eqref{eqn:appCovZterm2} becomes
\begin{align}
\label{eqn:appCovZterm4}
R_{ml,nk}&={\mathbb{E}}\left({\bs{v}}_{m}(i)\!+\!{\bs{v}}_{ml}^{(d)}(i)\!-\!{\bs{v}}_{ml,i}^{(u)}w^o\right)
\left({\bs{v}}_{m}(i)\!+\!{\bs{v}}_{mk}^{(d)}(i)\!-\!{\bs{v}}_{mk,i}^{(u)}w^o\right)^*
\left({\bs{u}}_{m,i}\!+\!{\bs{v}}_{ml,i}^{(u)}\right)^*\!\left({\bs{u}}_{m,i}\!+\!{\bs{v}}_{mk,i}^{(u)}\right)\nonumber\\
{}&={\mathbb{E}}\left({\bs{v}}_{m}(i)+{\bs{v}}_{ml}^{(d)}(i)\right)\left({\bs{v}}_{m}(i)+{\bs{v}}_{mk}^{(d)}(i)\right)^*
{\mathbb{E}}\left({\bs{u}}_{m,i}+{\bs{v}}_{ml,i}^{(u)}\right)^*\left({\bs{u}}_{m,i}+{\bs{v}}_{mk,i}^{(u)}\right)\nonumber\\
{}&\qquad\qquad+{\mathbb{E}}\,{\bs{v}}_{ml,i}^{(u)}w^ow^{o*}{\bs{v}}_{mk,i}^{(u)*}
\left({\bs{u}}_{m,i}\!+\!{\bs{v}}_{ml,i}^{(u)}\right)^*\!\left({\bs{u}}_{m,i}\!+\!{\bs{v}}_{mk,i}^{(u)}\right)\nonumber\\
{}&=\left(\sigma_{v,m}^{2}+\delta_{lk}\sigma_{v,ml}^{2}\right)\left(R_{u,m}+\delta_{lk}R_{v,ml}^{(u)}\right)
+\delta_{lk}w^{o*}R_{v,ml}^{(u)}w^oR_{u,m}+R_{v,ml}^{(u)}w^ow^{o*}R_{v,mk}^{(u)}\nonumber\\
{}&\qquad\qquad+\delta_{lk}\left({\mathbb{E}}\,{\bs{v}}_{ml,i}^{(u)*}{\bs{v}}_{ml,i}^{(u)}
w^ow^{o*}{\bs{v}}_{ml,i}^{(u)*}{\bs{v}}_{ml,i}^{(u)}\!-\!R_{v,ml}^{(u)}w^ow^{o*}R_{v,ml}^{(u)}\right)
\end{align}
where $\delta_{lk}$ denotes the Kronecker delta function. Evaluating the last term on RHS of \eqref{eqn:appCovZterm4} requires knowledge of the excess kurtosis of ${\bs{v}}_{ml,i}^{(u)}$, which is generally not available. In order to proceed, we invoke a separation principle to approximate it as
\begin{align}
\label{eqn:appCovZterm5}
{\mathbb{E}}\,{\bs{v}}_{ml,i}^{(u)*}{\bs{v}}_{ml,i}^{(u)}w^ow^{o*}{\bs{v}}_{ml,i}^{(u)*}{\bs{v}}_{ml,i}^{(u)}
\approx R_{v,ml}^{(u)}w^ow^{o*}R_{v,ml}^{(u)}
\end{align}
Substituting \eqref{eqn:appCovZterm5} into \eqref{eqn:appCovZterm4} leads to
\begin{align}
\label{eqn:appCovZterm6}
R_{ml,nk}&\!\approx\!\left(\sigma_{v,m}^{2}\!+\!\delta_{lk}\sigma_{v,ml}^{2}\right)\left(R_{u,m}\!+\!\delta_{lk}R_{v,ml}^{(u)}\right)
\!+\!\delta_{lk}\left(w^{o*}R_{v,ml}^{(u)}w^o\right)R_{u,m}\!+\!R_{v,ml}^{(u)}w^ow^{o*}R_{v,mk}^{(u)}\nonumber\\
&\!=\!\sigma_{v,m}^{2}R_{u,m}\!+\!R_{v,ml}^{(u)}w^ow^{o*}R_{v,mk}^{(u)}
\!+\!\delta_{lk}\left[(\sigma_{v,ml}^{2}\!+\!w^{o*}R_{v,ml}^{(u)}w^o)R_{u,m}\!+\!(\sigma_{v,m}^{2}\!+\!\sigma_{v,ml}^{2})R_{v,ml}^{(u)}\right]
\end{align}
From \eqref{eqn:appCovZterm3} and \eqref{eqn:appCovZterm6}, we get
\begin{align}
\label{eqn:appCovZterm7}
R_{ml,nk}&\approx R_{v,ml}^{(u)}w^ow^{o*}R_{v,mk}^{(u)}+\delta_{mn}\sigma_{v,m}^{2}R_{u,m}\nonumber\\
&\qquad\qquad+\delta_{mn}\delta_{lk}\!\left[(\sigma_{v,ml}^{2}\!+\!w^{o*}R_{v,ml}^{(u)}w^o)R_{u,m}
\!+\!(\sigma_{v,m}^{2}\!+\!\sigma_{v,ml}^{2})R_{v,ml}^{(u)}\right]
\end{align}
Substituting \eqref{eqn:appCovZterm7} into \eqref{eqn:appCovZterm1}, we obtain
\begin{align}
\label{eqn:appCovZterm8}
R_{z,lk}&\approx\left(\sum_{m\in\N_l}c_{ml}R_{v,ml}^{(u)}w^o\right)\left(\sum_{n\in\N_k}c_{nk}R_{v,mk}^{(u)}w^{o}\right)^*
+\sum_{m\in\N_l}\sum_{n\in\N_k}c_{ml}c_{nk}\delta_{mn}\sigma_{v,m}^{2}R_{u,m}\nonumber\\
{}&\qquad\qquad+\delta_{lk}\sum_{m\in\N_l}c_{ml}^2\left[(\sigma_{v,ml}^{2}+w^{o*}R_{v,ml}^{(u)}w^o)R_{u,m}
+(\sigma_{v,m}^{2}+\sigma_{v,ml}^{2})R_{v,ml}^{(u)}\right]
\end{align}
From \eqref{eqn:meanzkidef}--\eqref{eqn:zdef} and \eqref{eqn:bigSdef}--\eqref{eqn:Tkdef}, we arrive at expression \eqref{eqn:bigRzdef}.


\begin{thebibliography}{10}

\bibitem{Zhao12ICC}
X.~Zhao and A.~H. Sayed,
\newblock ``Combination weights for diffusion strategies with imperfect
  information exchange,''
\newblock in {\em Proc. IEEE Int. Conf. Commun. (ICC)}, Ottawa, Canada, June
  2012.

\bibitem{Tsitsiklis84TAC}
J.~Tsitsiklis and M.~Athans,
\newblock ``Convergence and asymptotic agreement in distributed decision
  problems,''
\newblock {\em {IEEE} Trans. Autom. Control}, vol. 29, no. 1, pp. 42--50, Jan.
  1984.

\bibitem{Tsitsiklis86TAC}
J.~Tsitsiklis, D.~Bertsekas, and M.~Athans,
\newblock ``Distributed asynchronous deterministic and stochastic gradient
  optimization algorithms,''
\newblock {\em {IEEE} Trans. Autom. Control}, vol. 31, no. 9, pp. 803--812,
  Sept. 1986.

\bibitem{Bertsekas97JOP}
D.~P. Bertsekas,
\newblock ``A new class of incremental gradient methods for least squares
  problems,''
\newblock {\em SIAM J. Optim.}, vol. 7, no. 4, pp. 913--926, 1997.

\bibitem{Nedic01JOP}
A.~Nedic and D.~P. Bertsekas,
\newblock ``Incremental subgradient methods for nondifferentiable
  optimization,''
\newblock {\em SIAM J. Optim.}, vol. 12, no. 1, pp. 109--138, 2001.

\bibitem{Rabbat05JSAC}
M.~G. Rabbat and R.~D. Nowak,
\newblock ``Quantized incremental algorithms for distributed optimization,''
\newblock {\em {IEEE} J. Sel. Areas Commun.}, vol. 23, no. 4, pp. 798--808,
  Apr. 2005.

\bibitem{Lopes07TSP}
C.~G. Lopes and A.~H. Sayed,
\newblock ``Incremental adaptive strategies over distributed networks,''
\newblock {\em {IEEE} Trans. Signal Process.}, vol. 48, no. 8, pp. 223--229,
  Aug. 2007.

\bibitem{Lopes06ASAP}
C.~G. Lopes and A.~H. Sayed,
\newblock ``Distributed processing over adaptive networks,''
\newblock in {\em Proc. Adapt. Sensor Array Process. Workshop (ASAP)}, MIT
  Lincoln Laboratory, MA, June 2006, pp. 1--5.

\bibitem{Lopes08TSP}
C.~G. Lopes and A.~H. Sayed,
\newblock ``Diffusion least-mean squares over adaptive networks: {F}ormulation
  and performance analysis,''
\newblock {\em {IEEE} Trans. Signal Process.}, vol. 56, no. 7, pp. 3122--3136,
  July 2008.

\bibitem{Cattivelli08ASILOMAR}
F.~S. Cattivelli and A.~H. Sayed,
\newblock ``Diffusion {LMS} algorithms with information exchange,''
\newblock in {\em Proc. Asilomar Conf. Signals, Systems, Computers}, Pacific
  Grove, CA, Oct. 2008, pp. 251--255.

\bibitem{Cattivelli10TSP}
F.~S. Cattivelli and A.~H. Sayed,
\newblock ``Diffusion {LMS} strategies for distributed estimation,''
\newblock {\em {IEEE} Trans. Signal Process.}, vol. 58, no. 3, pp. 1035--1048,
  Mar. 2010.

\bibitem{Cattivelli10TAC}
F.~S. Cattivelli and A.~H. Sayed,
\newblock ``Diffusion strategies for distributed {Kalman} filtering and
  smoothing,''
\newblock {\em {IEEE} Trans. Autom. Control}, vol. 55, no. 9, pp. 2069--2084,
  Sept. 2010.

\bibitem{Cattivelli08TSP}
F.~S. Cattivelli, C.~G. Lopes, and A.~H. Sayed,
\newblock ``Diffusion recursive least-squares for distributed estimation over
  adaptive networks,''
\newblock {\em {IEEE} Trans. Signal Process.}, vol. 56, no. 5, pp. 1865--1877,
  May 2008.

\bibitem{Li09SSP}
L.~Li and J.~A. Chambers,
\newblock ``Distributed adaptive estimation based on the {APA} algorithm over
  diffusion netowrks with changing topology,''
\newblock in {\em Proc. IEEE Workshop Stat. Signal Process. (SSP)}, Cardiff,
  UK, Aug. / Sept. 2009, pp. 757--760.

\bibitem{Takahashi10TSP}
N.~Takahashi, I.~Yamada, and A.~H. Sayed,
\newblock ``Diffusion least-mean squares with adaptive combiners: {F}ormulation
  and performance analysis,''
\newblock {\em {IEEE} Trans. Signal Process.}, vol. 58, no. 9, pp. 4795--4810,
  Sept. 2010.

\bibitem{Takahashi10ICASSP}
N.~Takahashi and I.~Yamada,
\newblock ``Link probability control for probabilistic diffusion least-mean
  squares over resource-constrained networks,''
\newblock in {\em Proc. IEEE Int. Conf. Acoust., Speech, Signal Process.
  (ICASSP)}, Dallas, TX, Mar. 2010, pp. 3518--3521.

\bibitem{Chouvardas11TSP}
S.~Chouvardas, K.~Slavakis, and S.~Theodoridis,
\newblock ``Adaptive robust distributed learning in diffusion sensor
  networks,''
\newblock {\em {IEEE} Trans. Signal Process.}, vol. 59, no. 10, pp. 4692--4707,
  Oct. 2011.

\bibitem{Tu11JSTSP}
S-Y. Tu and A.~H. Sayed,
\newblock ``Mobile adaptive networks,''
\newblock {\em {IEEE} J. Sel. Top. Signal Process.}, vol. 5, no. 4, pp.
  649--664, Aug. 2011.

\bibitem{Cattivelli11TSP}
F.~S. Cattivelli and A.~H. Sayed,
\newblock ``Modeling bird flight formations using diffusion adaptation,''
\newblock {\em {IEEE} Trans. Signal Process.}, vol. 59, no. 5, pp. 2038--2051,
  May 2011.

\bibitem{Li11EURASIP}
J.~Li and A.~H. Sayed,
\newblock ``Modeling bee swarming behavior through diffusion adaptation with
  asymmetric information sharing,''
\newblock {\em EURASIP J. Advances Signal Process.}, pp. 1--14,
  doi:10.1186/1687--6180--2012--18, Jan. 2012.

\bibitem{Towfic11MLSP}
Z.~Towfic, J.~Chen, and A.~H. Sayed,
\newblock ``Collaborative learning of mixture models using diffusion
  adaptation,''
\newblock in {\em Proc. IEEE Int. Workshop Machine Learn. Signal Process.
  (MLSP)}, Beijing, China, Sept. 2011, pp. 1--6.

\bibitem{Weng11Sensors}
Y.~Weng, W.~Xiao, and L.~Xie,
\newblock ``Diffusion-based {EM} algorithm for distributed estimation of
  {G}aussian mixtures in wireless sensor networks,''
\newblock {\em Sensors}, vol. 11, no. 6, pp. 6297--6316, June 2011.

\bibitem{Chen11TSP}
J.~Chen and A.~H. Sayed,
\newblock ``Diffusion adaptation strategies for distributed optimization and
  learning over networks,''
\newblock to appear in {\em IEEE Trans. Signal Process.}, 2012. [Also available online at http://arxiv.org/abs/1111.0034 as {\em arXiv:1111.0034v2 [math.OC]}, Oct. 2011.]

\bibitem{Jadbabaie03TAC}
A.~Jadbabaie, J.~Lin, and A.~S. Morse,
\newblock ``Coordination of groups of mobile autonomous agents using nearest
  neighbor rules,''
\newblock {\em {IEEE} Trans. Autom. Control}, vol. 48, no. 6, pp. 988--1001,
  June 2003.

\bibitem{Fax04TAC}
J.~A. Fax and R.~M. Murray,
\newblock ``Information flow and cooperative control of vehicle formations,''
\newblock {\em {IEEE} Trans. Autom. Control}, vol. 49, no. 9, pp. 1465--1476,
  Sept. 2004.

\bibitem{Olfati06TAC}
R.~Olfati-Saber,
\newblock ``Flocking for multi-agent dynamic systems: Algorithms and theory,''
\newblock {\em {IEEE} Trans. Autom. Control}, vol. 51, no. 3, pp. 401--420,
  Mar. 2006.

\bibitem{Olfati04TAC}
R.~Olfati-Saber and R.~M. Murray,
\newblock ``Consensus problems in networks of agents with switching topology
  and time-delays,''
\newblock {\em {IEEE} Trans. Autom. Control}, vol. 49, no. 9, pp. 1520--1533,
  Sept. 2004.

\bibitem{Barbarossa07SPM}
S.~Barbarossa and G.~Scutari,
\newblock ``Bio-inspired sensor netowrk design,''
\newblock {\em {IEEE} Signal Process. Mag.}, vol. 24, no. 3, pp. 26--35, May
  2007.

\bibitem{Nedic09TSP}
A.~Nedic and A.~Ozdaglar,
\newblock ``Distributed subgradient methods for multi-agent optimization,''
\newblock {\em {IEEE} Trans. Signal Process.}, vol. 54, no. 1, pp. 48--61, Jan.
  2009.

\bibitem{Dimakis10PROC}
A.~G. Dimakis, S.~Kar, J.~M.~F. Moura, M.~G. Rabbat, and A.~Scaglione,
\newblock ``Gossip algorithms for distributed signal processing,''
\newblock {\em Proc. {IEEE}}, vol. 98, no. 11, pp. 1847--1864, Nov. 2010.

\bibitem{Kar11TSP}
S.~Kar and J.~M.~F. Moura,
\newblock ``Convergence rate analysis of distributed gossip (linear parameter)
  estimation: Fundamental limits and tradeoffs,''
\newblock {\em {IEEE} J. Sel. Top. Signal Process.}, vol. 5, no. 4, pp.
  674--690, Aug. 2011.

\bibitem{Sayed08}
A.~H. Sayed,
\newblock {\em Adaptive Filters},
\newblock Wiley, NJ, 2008.

\bibitem{Abdolee11DCOSS}
R.~Abdolee and B.~Champagne,
\newblock ``Diffusion {LMS} algorithms for sensor networks over non-ideal
  inter-sensor wireless channels,''
\newblock in {\em Proc. IEEE Int. Conf. Dist. Comput. Sensor Systems (DCOSS)},
  Barcelona, Spain, June 2011, pp. 1--6.

\bibitem{Khalili11ACSP}
A.~Khalili, M.~A. Tinati, A.~Rastegarnia, and J.~A. Chambers,
\newblock ``Transient analysis of diffusion least-mean squares adaptive
  networks with noisy channels,''
\newblock {\em Int. J. Adapt. Control Signal Process.}, Sept. 2011,
  doi:~10.1002/acs.1279.

\bibitem{Khalili12TSP}
A.~Khalili, M.~A. Tinati, A.~Rastegarnia, and J.~A. Chambers,
\newblock ``Steady-state analysis of diffusion {LMS} adaptive networks with
  noisy links,''
\newblock {\em {IEEE} Trans. Signal Process.}, vol. 60, no. 2, pp. 974--979,
  Feb. 2012.

\bibitem{Tu11GlobeCom}
S-Y. Tu and A.~H. Sayed,
\newblock ``Adaptive networks with noisy links,''
\newblock in {\em Proc. IEEE Global Commun. Conf. (GLOBECOM)}, Houston, TX,
  Dec. 2011, pp. 1--5.

\bibitem{Kar09TSP}
S.~Kar and J.~M.~F. Moura,
\newblock ``Distributed consensus algorithms in sensor networks: Link failures
  and channel noise,''
\newblock {\em {IEEE} Trans. Signal Process.}, vol. 57, no. 1, pp. 355--369,
  Jan. 2009.

\bibitem{Mateos09EUROSIP}
G.~Mateos, I.~D. Schizas, and G.~B. Giannakis,
\newblock ``Performance analysis of the consensus-based distributed {LMS}
  algorithm,''
\newblock {\em EURASIP J. Advances Signal Process.}, pp. 1--19, Article ID
  981030, doi:10.1155/2009/981030, 2009.

\bibitem{Metropolis53JCP}
N.~Metropolis, A.~W. Rosenbluth, M.~N. Rosenbluth, A.~H. Teller, and E.~Teller,
\newblock ``Equations of state calculations by fast computing machines,''
\newblock {\em J. Chem. Phys.}, vol. 21, no. 6, pp. 1087--1092, 1953.

\bibitem{Xiao05IPSN}
L.~Xiao, S.~Boyd, and S.~Lall,
\newblock ``A scheme for robust distributed sensor fusion based on average
  consensus,''
\newblock in {\em Proc. ACM/IEEE Int. Conf. Inform. Process. Sensor Networks
  (IPSN)}, Los Angeles, CA, Apr. 2005, pp. 63--70.

\bibitem{Xiao04SCL}
L.~Xiao and S.~Boyd,
\newblock ``Fast linear iterations for distributed averaging,''
\newblock {\em System Control Lett.}, vol. 53, no. 9, pp. 65--78, Sept. 2004.

\bibitem{Jakovetic10TSP}
D.~Jakovetic, J.~Xavier, and J.~M.~F. Moura,
\newblock ``Weight optimization for consensus algorithms with correlated
  switching topology,''
\newblock {\em {IEEE} Trans. Signal Process.}, vol. 58, no. 7, pp. 3788--3801,
  July 2010.

\bibitem{Arenas05TIM}
J.~Arenas-Garcia, V.~Gomez-Verdejo, and A.~R. Figueiras-Vidal,
\newblock ``New algorithms for improved adaptive convex combination of {LMS}
  transversal filters,''
\newblock {\em {IEEE} Trans. Instrum. Meas.}, vol. 54, no. 6, pp. 2239--2249,
  Dec. 2005.

\bibitem{Mandic07ICASSP}
D.~Mandic, P.~Vayanos, C.~Boukis, B.~Jelfs, S.~L. Goh, T.~Gautama, and
  T.~Rutkowski,
\newblock ``Collaborative adaptive learning using hybrid filters,''
\newblock in {\em Proc. IEEE Int. Conf. Acoust., Speech, Signal Process.
  (ICASSP)}, Honolulu, HI, Apr. 2007, pp. 921--924.

\bibitem{Silva08TSP}
M.~T.~M. Silva and V.~H. Nascimento,
\newblock ``Improving the tracking capability of adaptive filters via convex
  combination,''
\newblock {\em {IEEE} Trans. Signal Process.}, vol. 56, no. 7, pp. 3137--3149,
  July 2008.

\bibitem{Kozat09ICASSP}
S.~S. Kozat and A.~C. Singer,
\newblock ``A performance-weighted mixture of {LMS} filters,''
\newblock in {\em Proc. IEEE Int. Conf. Acoust., Speech, Signal Process.
  (ICASSP)}, Taipei, Apr. 2009, pp. 3101--3104.

\bibitem{Candido10TSP}
R.~Candido, M.~T.~M. Silva, and V.~H. Nascimento,
\newblock ``Transient and steady-state analysis of the affine combination of
  two adaptive filters,''
\newblock {\em {IEEE} Trans. Signal Process.}, vol. 58, no. 8, pp. 4064--4078,
  Aug. 2010.

\bibitem{Boyd04}
S.~Boyd and L.~Vandenberghe,
\newblock {\em Convex Optimization},
\newblock Cambridge Univ. Press, 2004.

\bibitem{Horn85}
R.~A. Horn and C.~R. Johnson,
\newblock {\em Matrix Analysis},
\newblock Cambridge Univ. Press, Cambridge, UK, 1985.

\bibitem{Alnaffouri03TSP}
T.~Y. Al-Naffouri and A.~H. Sayed,
\newblock ``Transient analysis of data-normalized adaptive filters,''
\newblock {\em {IEEE} Trans. Signal Process.}, vol. 51, no. 3, pp. 639--652,
  Mar. 2003.

\bibitem{Laub05}
A.~J. Laub,
\newblock {\em Matrix Analysis for Scientists and Engineers},
\newblock SIAM, PA, 2005.

\bibitem{Tu11CAMSAP}
S-Y. Tu and A.~H. Sayed,
\newblock ``Optimal combination rules for adaptation and learning over
  netowrks,''
\newblock in {\em Proc. IEEE Int. Workshop Comput. Advances Multi-Sensor Adapt.
  Process. (CAMSAP)}, San Juan, Puerto Rico, Dec. 2011, pp. 317--320.

\end{thebibliography}
\end{document}